\documentclass[]{interact}

\usepackage{epstopdf}% To incorporate .eps illustrations using PDFLaTeX, etc.

\usepackage{natbib}% Citation support using natbib.sty
\bibpunct[, ]{(}{)}{;}{a}{}{,}% Citation support using natbib.sty
\theoremstyle{plain}% Theorem-like structures provided by amsthm.sty

\theoremstyle{definition}

\theoremstyle{remark}

\usepackage{subfig}
\usepackage{graphicx}
\DeclareMathOperator\arctanh{arctanh}
\def\bibfont{\small}%
\def\bibsep{\smallskipamount}%
\def\bibhang{24pt}%
\usepackage{placeins}
\begin{document}

%\articletype{ARTICLE TEMPLATE}

\title{Unifying warfighting functions in mathematical modelling: combat, manoeuvre, and C2}

\author{
\name{
Ryan Ahern\textsuperscript{a},
Mathew Zuparic\textsuperscript{b},
Keeley Hoek\textsuperscript{c} and
Alexander Kalloniatis\textsuperscript{b}
}
\affil{
\textsuperscript{a} 
Monash University, Clayton, VIC, Australia;
\textsuperscript{b} 
Defence Science and Technology Group, Canberra, ACT, Australia;
\textsuperscript{c}
Australian National University,
Canberra, Australia.
}
}
\maketitle
\begin{abstract}
The outcomes of warfare have rarely only been characterised by the quantity and quality of individual combatant force elements. The ability to manoeuvre and adapt across force elements through effective Command and Control (C2) can allow smaller or weaker forces to overcome an adversary with greater resource and fire-power. In this paper, we combine the classic Lanchester combat model with the Kuramoto-Sakaguchi model for phase oscillators on a network to create a flexible Networked-Lanchester-C2 representation of force-on-force military engagement. 
The mathematical model thus unifies three of the military
warfighting `functions': fires, manoeuvre and C2. We consider three illustrative use-cases, and show that an analytical treatment of a reduced model characterises global effects in the full system. For inhomogeneous forces we observe
that with appropriate balance between
internal organisational coupling, resource
manoeuvrability and even weaker lethality
the force can be adaptive to overcome
an initially stronger adversary. 
\end{abstract}


\maketitle
%%%%%%%%%%%%%%%%%%%%%%%%%%%%%%%%%%%%%%%%%%%%%%%%%%%%%%%%%%%%%%%%%%%%%%

\section{Introduction}
Since the dawn of organised warfare, strategy and tactics have enabled commanders to defeat numerically superior foes. In
particular, the ability to outmanoeuvre the forces of an opponent can open opportunities for a more decisive result than purely attritional warfare.
We propose to mathematically model these dynamics through coupled differential equations, integrating
combat, manoeuvre, and Command and
Control (C2).

%As a seminal example, \cite{v85} analysed the Napoleonic wars and highlighted the revolution to warfare brought about during these times by the introduction of a Headquarters, thus enabling greater levels of decentralised command, and superior manoeuvrability of the \textit{Grande Arm\'{e}e}:
%\begin{quote}
  %  Napoleon on campaign no longer attempted to keep the bulk of his forces concentrated under his own hand. No longer was the commander found doing everything important...Headquarters drastically cut down the burden of communicating and data processing that rested upon it. Since the corps were able to operate and hold out on their own, Imperial Headquarters could tolerate a far higher degree of uncertainty concerning the corps' momentary situation; this in turn made it possible to raise the level of performance to the point that it took the French armies only a few brief campaigns to overrun virtually an entire continent. 
    %
   % --- \cite{v85} p96--7
%\end{quote}
Napoleon is arguably the most famous, and defining, example of superior manoeuvrability of forces at a \textit{campaign} level
\cite{v85}, able to
 shift his forces between distributed and then concentrated at a decisive point, with synchronisation a key enabler for his victories.
There are also \textit{more localised} examples where a cleverly diffuse, agile and responsive force has been able to defeat a numerically superior opponent. 
At the battle of Cynoscephalae in 197 BC \cite{Polybius27}, the superior manoeuvrability and controllability of the \textit{maniple} tactical sub-division formation enabled Roman victory over the Macedonian \textit{phalanx} and heavy cavalry. 
At the battle of Cr\'{e}cy \cite{Reid07} in 1346, with the introduction of the English longbow, tactics were enabled to effectively disrupt the (usually devastating) heavy cavalry charge so that the English defeated a vastly numerically superior French army. 
Manoeuvre ideas were again
vindicated with the
Blitzkrieg of the German Wehrmacht in World War II, and coalition
forces in the First Gulf War.

%In this paper, we seek to differentiate these considerations of \textit{manoeuvrability}, \textit{disruption} and \textit{synchronisation} from physical war-fighting capability, and combine them in a coherent model which permits exploration of emergent effects.  

The model we propose builds upon the Lanchester square (directed-fire) model \cite{Lanchester},
\begin{eqnarray}
\dot{p}_B = -\alpha_{RB} p_R, \ \ \   p_B= p_{B0} \ \text{when}\ t=0 , \ \ 
\dot{p}_R = -\alpha_{BR} p_B, \ \ \   p_R= p_{R0} \ \text{when}\ t=0, 
\label{eq:Lanchester}
\end{eqnarray}
where $\alpha_{RB}$, $\alpha_{BR}$ give the combat effectiveness, or lethality, for the Red and Blue force, respectively, and $p_R$,$p_B$ are the resource \textit{populations} of Red and Blue respectively.
%Historically used to fit data for the battles of Ardennes \cite{Bracken95} and Kursk \cite{Lucas04}, the Lanchester square model approximates the simultaneous attrition of two opposing forces. 
%Though the previously cited authors, amongst others, have questioned the applicability of attrition in modelling contemporary engagements, we contend that the concept can still be well-founded in the modern-era, especially considering the shift towards robust distributed autonomous capabilities over single-entity assets. 
As a constant coefficient linear system of ordinary differential equations a conserved quantity from Eq.(\ref{eq:Lanchester}) is analytically solvable, meaning that lethalities and initial resources determine the outcome of the battle.
%
%via
%\begin{eqnarray}
%\begin{split}
%2 \left(\alpha_{BR} p_R \dot{p}_R -  \alpha_{RB} p_B \dot{p}_B \right) \equiv \frac{d}{dt}\left( \alpha_{BR} p^2_R - \alpha_{RB} p^2_B \right) =0,\\
%\Rightarrow  \alpha_{BR} p^2_R(t) - \alpha_{RB} p^2_B(t) = \alpha_{BR} p^2_{R}(0) - \alpha_{RB} p^2_{B}(0).
%\end{split}
%\label{conslaw}
%\end{eqnarray}

Note that other models 
lack continuous battle dynamics, as in the 
Colonel Blotto game \cite{Roberson06}, or the Hughes salvo model
\cite{Hughes95} of naval warfare.
The Lanchester model has been extended to fit empirical data in different historical contexts \cite{Bracken95} (see also \cite{Fricker98} and \cite{Lucas04}), has been extended to the more recent \textit{mixed forces} variant by  \cite{MacKay09} and \cite{Kress2018}, 
has been used for \textit{artificial intelligence} playing the real-time-strategy game \textit{Starcraft} \cite{Stanescu15} and for asymmetric guerilla warfare with  \cite{Deitchman62} and \cite{Schaffer68}.

In all these, the functions of command and manoeuvre are buried in the combat effectiveness {\it constant}. This is reasonable when fitting against data for two forces with
unchanged strategy over a drawn out engagement, but this is rarely true even for simple social animals,
as discussed in \cite{Adams03} and \cite{Plowes05}. It is also difficult to estimate {\it a priori} the combat effectiveness for a force, or assess the stability of the effectiveness to changes in organisation. 
More significantly, the role of information and communication technology in enabling the organisation of diverse (for example,
from different military services) elements of
a force has made the {\it network} paradigm
potent for contemporary and future warfare.
This has led, in recent decades, to
the proposition of what was
coined as `Network Centric Warfare'
\cite{Alberts1999} 
linking network connectivity,
self-synchronisation of force elements
and effective military outcomes. 
Though the Lanchester model has been
treated as an abstraction for
network {\it cyber} attack and defence
\cite{LiuEtAl2013}, a model with explicit network structure
was only proposed in \cite{Kall19} by generalising
MacKay's  \cite{MacKay09} 
mixed-forces model. Here network optimisation led to identification of concepts
in the Manoeuvre Theory of Warfare.

In this paper we bring C2 into
the formulation of \cite{Kall19} using
the Kuramoto-Sakaguchi
\cite{KuramotoSakaguchi} model
of phase oscillators (which builds on \cite{Kuramoto84}),  thereby unifying
three key warfighting functions, C2, Manoeuvre and
Fires in a single mathematical model.
For general networks,
as developed by 
 \cite{Acebron05}, \cite{Arenas08}, \cite{Dorfler14} and \cite{Rodrigues16} 
 the model is
\begin{equation} \label{eq:Kuramoto}
\dot{\theta_i} = \omega_i - \sum_{j \in {\cal K}} \sigma_{ij}  {\mathcal{K}_{ij}} \sin(\theta_i-\theta_j-\Phi_{i j}), \;\; i \in \mathcal{K},
\end{equation}
where for the $i$th node, $\theta_i \in \mathbb{S}^1$ is the phase, $\omega_i \in \mathbb{R}$ is the natural frequency, $\sigma_{ij}\in \mathbb{R}_+$ is the coupling strength on node $i$ from node $j$, $\cal{K}$ is the adjacency matrix representing the network, and $\Phi_{ij}\in \mathbb{S}^1$ represents a `frustration' between node $i$ and $j$, where phase $\theta_i$ attempts to stay amount $\Phi_{ij}$ ahead of phase $\theta_j$. 
It is well-known that at sufficiently
high couplings $\sigma_{ij}$, phases approximately
align $\theta_i(t)\approx \theta_j(t)$ and 
synchronise to (for undirected graphs)
the mean of the natural frequencies $\bar{\omega}$.

For application to C2
\cite{Kalloniatis-McLenn-Rob2019}, the 
phases represent decision
cycles of individual agents,
the natural frequencies the decision speeds
of agents if left to themselves,
coupling strengths represent the tightness
of relationship (or responsiveness to
changes in decision state) between agents,
the network represents the formal
and informal organisational C2 structure, and the frustrations
represent how far ahead 
agents {\it seek} to be in relation to
other agents. The same
framework of Eq.(\ref{eq:Kuramoto})
can be extended to two (or more \cite{ZuparicEtAl2019}) forces or populations
in competition with each other.
Two elements quantified
here mathematically are well documented
as {\it qualitative} properties of
decision-making and organisations:
the cyclicity of individual decision-making, in the
Perception-Action cycle of
cognitive psychology \cite{Neisser76} or the
Observe-Orient-Decide-Act (OODA) loop of military and business strategy
\cite{Osinga06}; and
the role of loose and tight coupling
in organisations as articulated by \cite{Weick1976}, \cite{Perrow2011} and \cite{HollenSpitz2012}.
This model thus brings
long articulated ideas in organisational theory into
a dynamical mathematical model.

To illustrate these ideas more specifically
in terms of the variables we have introduced
consider that the first $N$ components of
$\theta_i$ represent the decision state in
a continuous version of the OODA loop for
individual Blue members of the Blue
C2 system or network. The remaining $M$ components
represent the (OODA) decision state of Red agents
within their own corresponding C2 system. In the
spirit of Boyd's admonition each of Blue and Red
will seek to stabilise their decision state
ahead of each other, $\phi_i > \phi_j$ 
where $i,j$ may represent Blue or Red. 
Thus each side aims for locking a certain
phase shift ahead of the other: $\Phi_{BR}$ for Blue
in relation to Red, and $\Phi_{RB}$ for Red
with respect to Blue; for simplicity here
we treat the frustration values homogeneously
within each group, Blue or Red. These frustration
values represent a target state for each side.
Whether they achieve that depends on the dynamics
represented in Eq.~\ref{eq:Kuramoto}, subject
to their couplings, their connectivity, and their
individual capacity for rapid decision making.

{Till now, the `decisions' referred to have been in
the abstract. In this paper} we
wed this to what C2 is meant to accomplish
in a military force. Thus we model
the situation that if one of Blue or Red have
their phase $\theta_i$ ahead of the other,
$\theta_j$ they will gain an advantage in the
Lanchester dynamics. We thus propose a mathematical
model naturally unifying the combat, manoeuvre and C2 warfighting functions using a multi-layer network formulation \cite{MultiNetwork}. The model is more compact than an alternative contemporary equation-based
approach by \cite{McLemore2016}, noting that they spatially embed force elements and omit C2.

For homogeneous forces our unified model can be simplified to a global form, with which we start to illustrate the construction, and build
on a use-case from experience with the Kuramoto-Sakaguchi model \cite{PhysicaA}. Though not solvable,
the global model permits considerable dimensional reduction and insight. 
When inhomogeneity is high the global model is more limited, but the full model demonstrates {\it adaptive} dynamics where one force, say Blue, must change its resource distribution to adjust for a stronger Red force. Finally, with a fictitious complex network we illustrate how
the complete model may be used in a classic
Operational Research sense, to investigate trade-offs between investment in 
fire-power or C2 systems.

In the following, we present the principles of the model by showing how
C2 and combat dynamics may be combined
in a unified Lanchester-Kuramoto-Sakaguchi model.
We examine this model only insofar as it
builds an intuition and approximation
scheme that can be applied to the
full model. We then formulate the
fully networked Lanchester-Kuramoto-Sakaguchi model
and consider it for the first two use-cases, followed by
a complex network to show the utility of the
networked model in Operational Research.
We conclude and discuss future 
developments. Further details of analytical calculations
and examination of a use-case are relegated to the Appendix.

\section{Unifying the Lanchester equations with the Kuramoto-Sakaguchi model}
\label{sec2}
In the following we first outline the adversarial form 
of the Kuramoto-Sakaguchi model before then showing it may be
unified with the Lanchester equations and then
analysing to some degree the properties of the model.
We emphasise that this is both to establish the 
principles of how we develop the more advanced model
later on, and establish some intuition into its behaviours
to provide a foundation for analysing the fully networked
model in Sec.3. We provide two cases of such
a model, one where the C2 impacts on the Lanchester 
dynamics unilaterally, and the other where
attrition can in turn degrade the ability of a force
to undertake its decision making. 
We examine how far a semi-analytical approach can
go in understanding behaviours detected in numerical
solutions.

%This formulation causes node $i$ to attempt to deviate from its natural frequency to match the phase target of connected nodes, which in the case of a homogeneous network coupling $\sigma_{ij}=\text{const} \in \mathbb{R}\ \forall\ i,j \in \text{K}$ is a mutual effect between all connected nodes. For example, in the case of two nodes,
%$$
%\text{For } \dot{\theta_i} \
%\begin{cases}
%\text{If}\ \theta_i > \theta_j+\Phi_{ij}, &\ -\text{sgn}[\sigma_{j i}  \text{M}_{i j} \sin(\theta_i-\theta_j-\Phi_{i j}}] < 0 \quad \therefore \dot{\theta_i} < \omega_i \\
%\text{If}\ \theta_i = \theta_j+\Phi_{ij}, &\ -\text{sgn}[\sigma_{j i}  \text{M}_{i j} \sin(\theta_i-\theta_j-\Phi_{i j}}] = 0 \quad \therefore \dot{\theta_i} = \omega_i \\
%\text{If}\ \theta_i < \theta_j+\Phi_{ij}, &\ -\text{sgn}[\sigma_{j i}  \text{M}_{i j} \sin(\theta_i-\theta_j-\Phi_{i j}}] > 0 \quad \therefore \dot{\theta_i} > \omega_i 
%\end{cases}
%$$
\subsection{Adversarial form of Kuramoto-Sakaguchi model}
To begin with we show how Eq.(\ref{eq:Kuramoto}) may be written
to distinguish two adversary forces
that we call the `Blue-Red' formulation \cite{PhysicaA}. 
We adopt here a more
compact formulation compared to that work with
the following block matrix variable form for coupling,
frustration and network
structures:
\begin{equation}
\sigma = \left[ 
\begin{array}{cc}
\sigma_B & \zeta_{BR}   \\ 
\zeta_{RB} & \sigma_R  
\end{array} \right], \;\;
\Phi = \left[ 
\begin{array}{ccc}
0 & \phi^{BR} \\ 
\phi^{RB} & 0  
\end{array} \right], \;\;
{\cal K} = \left[ 
\begin{array}{cc}
\mathcal{B} & \mathcal{A}^{BR}  \\ 
\mathcal{A}^{RB} & \mathcal{R}
\end{array} \right]. 
\label{Kuradef}
\end{equation}
The blocks in $\sigma$ and $\Phi$ have the same
dimension as those in ${\cal K}$ so that diagonal unit matrices
are implicit in the respective blocks, while
$\theta$ becomes a vector composed of two blocks
of dimension $|{\cal B}|$ and $|{\cal R}|$ for the number
of Blue and Red agents, respectively.
Here $\zeta$ denotes an inter-network coupling, $\mathcal{A}$ matrices encode any inter-network edges and the frustration $\phi$ is set to zero for intra-network connections --- though this condition can be relaxed in future. 
%to explore the effects of intra-network competition. 
In this notation $\zeta_{BR}$ represents
Red's coupling to Blue, $\phi^{BR}$ the
degree to which Red seeks to be ahead of Blue's phase and $\mathcal{A}^{BR}$ the structure of
Red's interactions with Blue; reversing
the labels $B$ and $R$ gives Blue's
pattern of engagement with Red.
We thus capture a distinction
between internal processes of a force, Blue
or Red, from their interactions with each other. In particular,
the vanishing diagonals, and
taking $\{\phi^{BR},\phi^{RB}\}>0$
means that Blue and Red seek to internally
synchronise their decision-cycles while
trying to be ahead of each others cycles.
Such multi-network models \cite{MultiNetwork} exhibit collective synchronisation  effects that motivate their use in a C2-enabled Lanchester model.

In this paper, local synchronisation denotes the extent to which all nodes of a particular sub-network are phase-locked, i.e. $\theta_i=\theta_j$ $\forall$ $i,j \in \mathcal{B}$ or ${\cal R}$ describes a phase-locked Blue or Red sub-network. 
(If $\theta_i=\theta_j$ for all nodes across the combined Blue-Red system then we may speak of global synchronisation.)
We may quantify this by sampling the instantaneous agreement in phase, referred to as the `order parameter'
\begin{equation}
O_{B}(t) \equiv \frac{\left|\sum_{j\in\mathcal{B}}e^{\sqrt{-1}\theta_j(t)} \right|}{|{\cal B}|}, \;\; O_{R}(t) \equiv \frac{\left|\sum_{j\in\mathcal{R}}e^{\sqrt{-1}\theta_j(t)} \right|}{|{\cal R}|}  . \label{eq:OrderParam}
\end{equation}
Note here we explicitly use the imaginary unit
$\sqrt{-1}$ rather than the oft-used symbol
$i$ because of the prevalence of the latter 
as an index for graph structure. Also, factors such as
$\sum_{j\in\mathcal{B}}e^{\sqrt{-1}\theta_j(t)}$
are easily seen to amount to the cosine of
difference of phases through multiplication 
by the complex conjugate and symmetry considerations.
Thus, the range of $O$, for either Blue or Red, is 
seen to be between zero and unity inclusive, where $O=1$ is attained if the phase of all nodes coincide 
(zero phase difference) and $O=0$ if the node phases are randomly distributed uniformly around the decision circle. Thus, the order parameter allows us to represent force cohesion and include this effect in the organisational component of force effectiveness.

\subsection{Unifying Lanchester and Kuramoto: C2 as Force Multiplier}
The key to unifying the two dynamical modelling paradigms is the notion that
C2 can be a force-{\it multiplier}, namely 
a factor that makes a set of weapons
act as more than just the sum of their
individual strengths. 
Thus, in this paper, we assume that force effectiveness in Lanchester is 
{\it multiplicatively separable} into  \textit{physical} $\kappa$ and \textit{organisational} $\Omega$ components, and apply the Kuramoto-Sakaguchi \cite{KuramotoSakaguchi} model as a framework to introduce organisational decision-making type dynamics. 
Hence, each lethality coefficient is 
factored as
 \begin{equation}
\alpha_{BR} = \kappa_{BR}\cdot\Omega_{BR}, \quad \alpha_{RB} = \kappa_{RB}\cdot\Omega_{RB}. \label{eq:SeperableCombatCoef}
\end{equation}
For the time being, we keep the physical effects $\kappa \in \mathbb{R}_+$ as a constant positive variable, but generate the organisational effects $\Omega$ dynamically, $\Omega=\Omega(t)$. These effects are achieved via the Red-Blue Kuramoto-Sakaguchi model originally considered in \cite{PhysicaA}
 \begin{eqnarray} 
\begin{split}
\dot{\theta}_{i}^B = \omega_i^B - \sigma_B \sum_{j \in \mathcal{B}} \mathcal{B}_{i j}\sin(\theta_i^B-\theta_j^B) - \zeta_{BR} \sum_{j \in \mathcal{R}} \mathcal{A}_{i j}^{BR}\sin(\theta^B_i-\theta^R_j-\phi^{BR}),\;\; i \in {\cal B},\\
\dot{\theta}_{i}^R = \omega_i^R - \sigma_R \sum_{j \in \mathcal{R}} \mathcal{R}_{i j}\sin(\theta_i^R-\theta_j^R) - \zeta_{RB} \sum_{j \in \mathcal{B}} \mathcal{A}_{i j}^{RB}\sin(\theta^R_i-\theta^B_j-\phi^{RB}), \;\; i \in {\cal R}
\end{split}\label{eq:BlueKuramoto}
\end{eqnarray}
where ${\theta}_{i}^B$ and ${\theta}_{i}^R$ are time-dependent. 
In the following we seek a form for $\Omega$ that captures the advantage provided by
a force being collectively ahead of the decision cycle
of the adversary.

%A first step to such a model can be one where
%the Blue and Red forces only conduct their
%internal C2 through the pure Kuramoto 
%equations
%and the impact of organisational
%effectiveness in the Lanchester model
%is through the order parameters for
%each force
%$\Omega_B=O_B, \Omega_R=O_R$.
%This compactly represents a detached
%headquarters enabling coherence of fires
%of a homogeneous force, but captures none of
%the dynamics of decision-advantage
%in the original OODA model.

\subsection{Global phases and organisational effect on combat}
\begin{figure}[h]
\centering
\subfloat[][Decision advantage on $\mathbb{S}^1$]{\includegraphics[width=0.4 \linewidth]{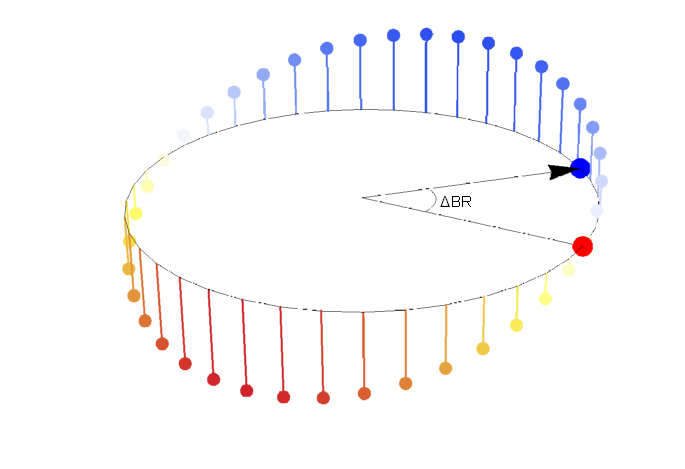}}
~
\subfloat[][Potential decision advantage functions ]{\includegraphics[width=0.6\linewidth]{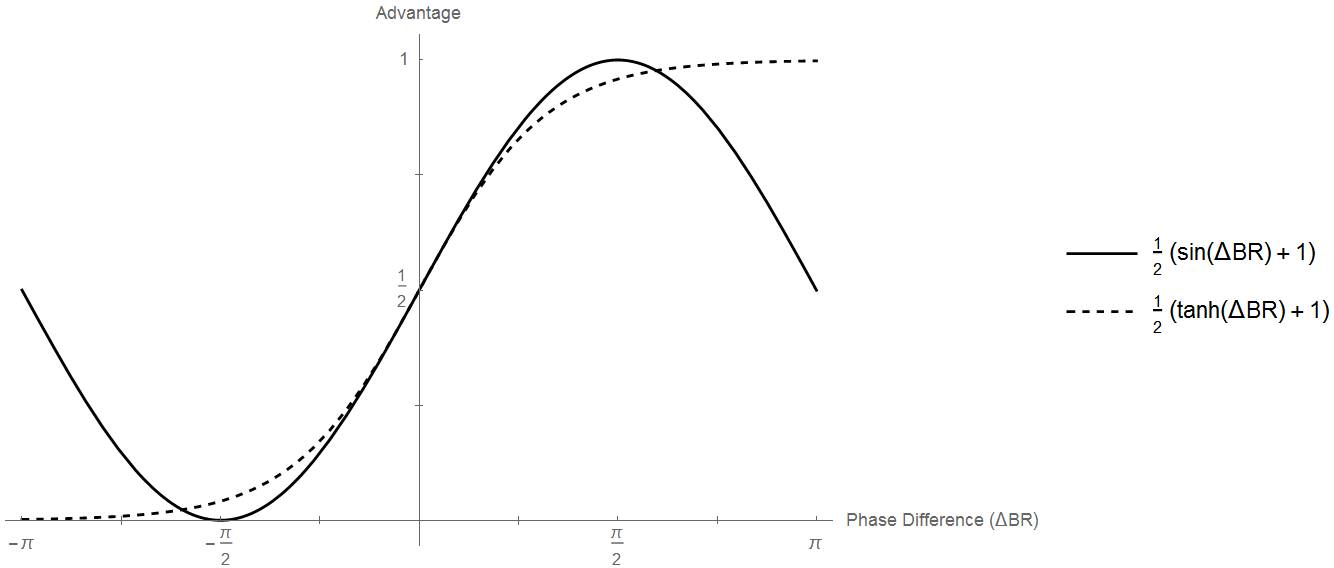}}
\caption{A three-dimensional perspective of the interpretation of sinusoidal advantage location on the decision circle (a) and a plot of different advantage functions (b). In both cases the Red node/cluster is located at zero angle and the Blue node/cluster is swept through values of Phase Difference $\Delta_{BR}$.}
\label{fig:PiMPiAdv}
\end{figure}

Because the Lanchester model in Eq.(\ref{eq:Lanchester}) only describes homogeneous forces,
we need to extract correspondingly global properties
of the adversarial C2 model. To this end we consider the \textit{global phases} (or centroids)
for agents in the Blue and Red networks, denoted by $\Theta_B$ and $\Theta_R$, respectively,
\begin{equation}
\Theta_B = \frac{1}{|{\cal{B}}|} \sum_{i \in {\cal{B}} }\theta^B_i, \;\; \Theta_R = \frac{1}{|{\cal{R}}|} \sum_{j \in {\cal{R}} }\theta^R_j.
\end{equation}
The difference between global phases
is designated $\Delta_{BR}$:
\begin{equation}
\Delta_{BR} \equiv \Theta_B - \Theta_R = - \Delta_{RB}.
\end{equation}
Applying the difference of global phases, we consider the following form for the organisational component of Blue/Red force effectiveness against Red/Blue,
 \begin{equation} 
 \Omega_{BR} = O_{B}\left( \frac{1+ \sin\Delta_{BR} }{2} \right) , \;\; \Omega_{RB} = O_{R}\left( \frac{1+ \sin\Delta_{RB} }{2} \right).
\label{eq:OrgComponent}
 \end{equation}
These $\Omega$ 
%in Eq.(\ref{eq:OrgComponent}) 
consist of the multiplication of two dynamic quantities 
%--- $O$ and $(1+\sin \Delta)/2$ --- 
which vary between zero and unity. Hence, each side gains organisational advantage by maintaining good intra-network synchronisation, \textit{and} staying ahead (maximally $\pi/2$) of their competitor's collective decision phase. Conversely, organisational advantage is lost by the loss of intra-network synchronisation, \textit{and/or} being behind (with a minimum achieved at $-\pi/2$) their competitor's collective decision phase, as depicted in Fig.\ref{fig:PiMPiAdv}. Additionally, the choice of sine (as opposed to tanh) allows decision advantage to be periodic and thus continuous over domain $\mathbb{S}^1$.
%, whereas tanh creates a discontinuous relationship which is both unphysical and numerically difficult to implement.
 %
Inserting this into Eq.(\ref{eq:SeperableCombatCoef}) yields
%the following set of Lanchester expressions,
 \begin{equation}
 \dot{p}_B = -\kappa_{RB}O_{R}\left( \frac{1- \sin \Delta_{BR} }{2}\right) p_R \mathcal{H}(p_B), \;\;  \dot{p}_R = -\kappa_{BR}O_{B}\left( \frac{1+ \sin \Delta_{BR} }{2}\right) p_B \mathcal{H}(p_R),
\label{fullLan}
 \end{equation}
where in order to avoid the appearance of negative populations, we enforce non-negativity with a Heaviside-step function $\mathcal{H}(x)$. The full generalised C2-Lanchester model is thus Eq.(\ref{fullLan}), in combination with Eq.(\ref{eq:BlueKuramoto}) representing the
dynamics of adversarial organisational decision cycles. We provide a comprehensive list of variables and parameters used to define the Lanchester-Kuramoto-Sakaguchi model in Table \ref{tab:tab1}.
%\begin{equation}
%\mathcal{H}(x) \equiv \begin{cases}
%1, &\ x>0\\
%0, &\ x \leq 0
%\end{cases}.
% \label{Heavi}
%\end{equation}
% Thus, the organisational component allows a force to access the maximum possible physical combat effectiveness ($\alpha = \kappa$) if it is both cohesive ($O=1$) and ahead of the opposition ($\Delta = \pi/2$). 

\begin{table}[ht]
\caption{Summary of the Kuramoto-Sakaguchi and Lanchester inspired variables and parameters defined in Sec. \ref{sec2} and \ref{sec3}, and their physical interpretations.} % title of Table
\centering % used for centering table
\begin{tabular}{c c c } % centered columns (4 columns)
\hline %inserts horizontal line
expression & name & interpretation \\ [0.5ex] % inserts table
%heading
\hline % inserts single horizontal line
$\{\theta^B,\theta^R\}$ & phase & agent decision state   \\
$\{{\cal B}, {\cal R}\}$ & adjacency matrix & internal decision-making network \\
$\omega^X$ & natural frequency &  decision-speed of ${\cal X}$ agent in isolation  \\
$\sigma_X$ & intra-network coupling & intensity of intra-agent interaction in ${\cal X}$\\
$ {\cal A}^{XY}$ & inter-adjacency matrix & topology between ${\cal X}$ and ${\cal Y}$ decision-makers   \\
$ \zeta_{XY}$ & inter-network coupling & intensity of ${\cal X}$'s interaction with ${\cal Y}$ \\ 
$ \phi^{XY}$ & frustration & ${\cal X}$'s strategy against ${\cal Y}$'s decision-makers \\ 
$O_X$ & global order parameter & measure of global synchronisation of\\
&& decision-making agents in ${\cal X}$ \\
$O^X_i$ & local order parameter & measure of local synchronisation of agent $i$\\
&& with its nearest neighbours in ${\cal X}$ \\
$\Theta_X$ & global phase & average phase of agents in ${\cal X}$ \\
$\Delta_{XY}$ & average phase difference & average difference of phases between\\ && networks ${\cal X}$ and ${\cal Y}$ \\
$\{p^B, p^R \}$ & population size & force/resource strength \\
$\kappa_{XY}$ & physical effectiveness & efficacy of ${\cal X}$'s physical\\
&& capabilities applied to $\cal Y$\\
$\Omega_{XY}$ & organisational effectiveness & efficacy of ${\cal X}$'s decision-making\\
&&capabilities applied to $\cal Y$ \\
$\alpha_{XY}$ & lethality & ${\cal X}$'s total combat effectiveness against ${\cal Y}$\\
$\{{\cal M}^B, {\cal M}^R\}$ & manoeuvre networks & internal force-flow network topology\\
${\cal E}^{XY}$ & engagement network & topology between ${\cal M}^X$ and ${\cal M}^Y$\\
&& force-networks where combat occurs \\
$\delta$ & flow moderation & dynamic moderation of force transfer\\
d & engagement moderation & dynamic moderation of combat outcomes\\
[1ex] % [1ex] adds vertical space
\hline %inserts single line
\end{tabular}
\label{tab:tab1} % is used to refer this table in the text
\end{table}

\subsection{Numerical simulation: use-case 1} \label{sec:NSofPhysA}
We now turn to numerical simulation of the model
to gain insight into its behaviours.
For this purpose
we approximate the Heaviside function as follows:
\begin{equation}
\mathcal{H}(x)
\ \approx
\left[ 1+ \tanh\left(
(x-\epsilon_1)/\epsilon_2\right) \right]/2
% \frac{1+\tanh %\frac{x-\epsilon_1}{\epsilon_2} %}{2},
\label{Heavi2}
\end{equation}
where $\epsilon_1 > \epsilon_2 > 0$ are small. 
In this paper we numerically solve using Mathematica's NDSolve method, to which end we choose $\epsilon_1=10^{-15}, \epsilon_2=10^{-20} $. 
Here we consider a use-case of competition between agents
organised, respectively, in a hierarchical and a random 
organisational structures, as in \cite{Holder17}
where we have already established a thorough
understanding of the potential dynamics of the model. 
These two organisations may be seen as caricatures
of the formal chain-of-command of traditional military forces, on
the one-hand, and the {\it ad hoc} nature of, say, many insurgent
groups. This simplification is
purely for illustrative purposes.

%For an in-depth discussion on the synchronisation, stability and fragmentation characteristics of the Kuramoto-Sakaguchi model in Eq.(\ref{eq:BlueKuramoto}), readers are directed to that paper. We shall briefly introduce the scenario and relevant parameters. 

\begin{figure}[h]
\centering
\subfloat[][Blue Sub-network]{\includegraphics[width=0.45\linewidth]{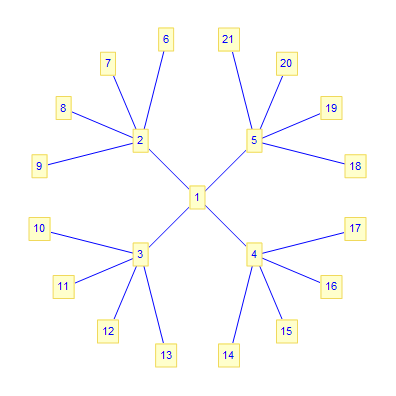}}
~
\subfloat[][Red Sub-network]{\includegraphics[width=0.45\linewidth]{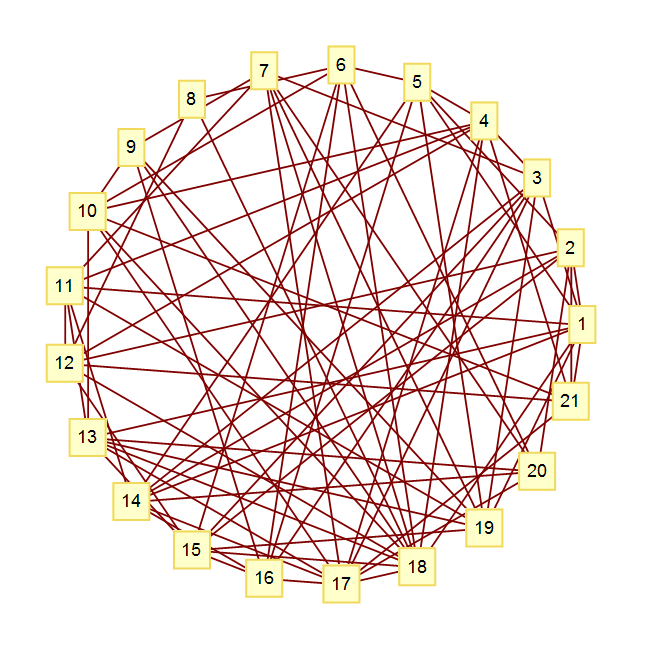}}
\caption{A graphical depiction of the sub-network structure for each force in use-case 1. The first five nodes, labelled 1-5, have no inter-network connections. All other nodes, labelled 6-21, are connected to the identically labelled node in the opposing network.}
\label{fig:NetworksPlot}
\end{figure}

Consider the nearly symmetric
scenario of \cite{PhysicaA}
\begin{equation}
|{\cal B}| = |{\cal R}|,\;\; \zeta_{BR} = \zeta_{RB},
\end{equation}
yet with different internal networks, couplings, frustrations and physical combat effectiveness. The sub-network structures themselves take the form described in Fig.\ref{fig:NetworksPlot} with Blue a hierarchy defined by a complete 4-ary tree such that $N_B=21$ and Red a random Erd\"{o}s-R\'{e}nyi network with an edge probability of 0.4 between the same number $N_R=21$ nodes. 
For the adversarial engagement between Blue and Red, we take that one each of the leaf nodes of
Blue (the tree, labelled $6,\dots,21$) engages with a correspondingly indexed
member of Red.

We choose internal couplings $\sigma_B = 8,\ \sigma_R=0.5$ sufficiently large that each network achieves internal synchronisation when considered separately. 
Again, this crudely represents the property that 
traditional military hierarchies also apply strong
discipline and training to achieve tight interaction
between members, while {\it ad hoc} organisations
are much `looser' affairs. That such
couplings are necessary for the respective structures
to achieve internal coherence is
confirmation that the Kuramoto model captures aspects
of real world C2 dynamics.  Note that, even though
this appears to negate the difference between
the Blue and Red we shall see that the higher
connectivity of the latter leaves its imprint in
many of the results.

\begin{figure}[h]
\centering
\includegraphics[width=1\linewidth]{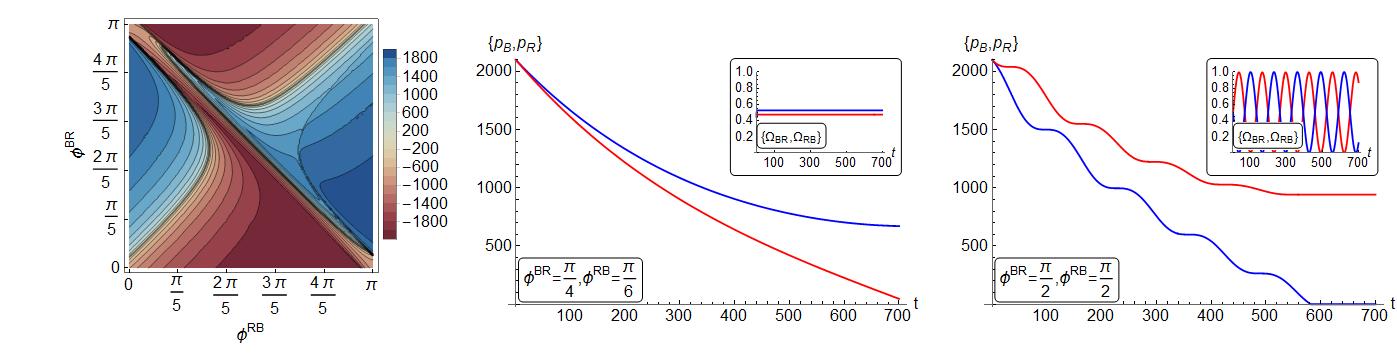}
\caption{Left most panel presents contour plot for the value of $p_{final}$ generated from  Eqs.(\ref{eq:BlueKuramoto},\ref{fullLan}), varying the frustration parameters $\phi^{BR}$ and $\phi^{RB}$. 
Colours are chosen to match where Blue or Red win the engagement. The middle and right most panels give examples of the the Lanchester trajectories for specific choices of frustration parameters. Additionally, the insets of the middle and right most plots show the dynamics with respect to time of the organisational component of each of Blue and Red's effectiveness, given in Eq.(\ref{eq:OrgComponent}).}
\label{fig:CombinationPlot}
\end{figure}

We also set $\zeta_{BR}=\zeta_{RB}=0.4,\ \kappa_{BR}=\kappa_{RB}=0.005,\ $ and randomly select the natural frequencies from a uniform distribution $[0, 1]$. The
particular choice we have worked
with in \cite{PhysicaA} has $\bar{\omega}_R = 0.551,\ \bar{\omega}_B = 0.503$, where $\bar{\omega}_B$ and $\bar{\omega}_R$ are the mean of the Blue and Red natural frequencies respectively. This small difference in means means
that in the absence of inter-network interactions Red would slowly lap Blue. 
Typically, random initial conditions are chosen for
numerical calculations which typically thermalise
in pure Kuramoto dynamics. However, in view of
the role that initial conditions play in the Lanchester
model, as in its conservation laws,
we choose an unconventional set of initial conditions here
to minimise {\it random} advantages to one side or the other.
We choose the phases of all nodes to be initially equi-spaced in the circular sector enclosed by $[-\frac{\pi}{4},\frac{\pi}{4}]$, i.e. $\theta_{0, j}^B = \theta_{0, j}^R = -\frac{\pi}{4}+\frac{\pi}{2}\frac{j-1}{|{\cal B}|-1}$ where $i = 1, 2, ..., |{\cal B}|$. The starting population for each side is equal $\left. p_B \right|_{t=0} = p_{B0} = \left. p_R \right|_{t=0}  = p_{R0} = 2100$. The system is solved numerically by Mathematica's NDSolve until $t_{final} = 10^4$, at which point the difference combat populations
\begin{equation}
 p_{final} = \left. p_B \right|_{t=t_{final}} - \left. p_R \right|_{t=t_{final}}
 \label{pfinal1}
\end{equation}
 is recorded, as shown in Fig. \ref{fig:CombinationPlot}.

The left hand panel of Fig.\ref{fig:CombinationPlot} gives the contour plot for $p_{final}$ of value for value of $ p_{final}$ generated from  Eqs.(\ref{eq:BlueKuramoto},\ref{fullLan}), while varying the frustration parameters $\phi^{BR}$ and $\phi^{RB}$.
We colour code to highlight where one side or the other
wins the engagement. We note although the contour plot does possess some symmetry, Red generally does better than Blue
by possessing a wider range of parameter across the
space of values that provide for Red victory. 
This is a consequence of Red's more connected
network structure creating a degree of 
greater robustness even though the couplings
of Blue and Red, respectively, are such as to
give the more poorly connected Blue a `fighting chance',
reflected in the degree of symmetry seen in
the left-most panel of Fig.\ref{fig:CombinationPlot}. The majority of the asymmetry is focused in the diagonal connecting the points $(0, \pi)$ and $(\pi,0)$; along this diagonal Red does comparatively better than Blue. By visual inspection we observe that along this diagonal the \textit{cross over} from favourable outcomes for Red to Blue is very sharp, indicated by the bold black lines. These sharp cross overs are in contrast to the gradual contour changes in other regions of the map. 

In the middle and right panels we present Lanchester trajectories for two specific choices of frustration parameters. The middle panel's trajectories, with choice $(\phi^{BR}, \phi^{RB}) = (\pi/4, \pi/6)$, decay approximately exponentially with
Blue winning the engagement. The trajectories in the right most panel, for $(\phi^{BR}, \phi^{RB}) = (\pi/2, \pi/2)$, show quasi-periodic decay but with Red victory. The insets in both these panels show the value of the organisational effects over time reflected in the $\Omega$ factors that
themselves depend on the C2 dynamics. Indeed,
the slightly better and flat synchronisation of Blue in
the middle panel is consistent with the decay of the
Lanchester trajectory. Similarly, the oscillations in the
right-most panel are consistent with those
in the C2 dynamics in the inset where both sides achieve internal
synchronisation, but their global phases periodically oscillate. There is no locking between Blue and Red with one or the other sustaining a lead in decision
state. In other words, as Blue and Red
agents pass each other in decision state they
momentarily gain an attritional advantage seen in
the oscillations in $p_B,p_R$. The fact that Red 
nonetheless wins that engagement
needs deeper explanation, to which we turn now.

\subsection{Clustering assumption and simplification}
Given that the model detailed in Eqs.(\ref{eq:BlueKuramoto},\ref{fullLan}) is a coupled system of nonlinear ordinary differential equations, generally it may only be solved numerically. Nevertheless, as the form of the Kuramoto model is currently unaffected by the Lanchester trajectories, we may leverage the work done on understanding the stability and synchronisation properties of two-network adversarial models with frustrations \cite{PhysicaA}. Beginning with the assumption that the Blue and Red networks are approximately phase synchronised via
\begin{equation}
\theta^B_i \approx b_i + \Theta_B, \;\; \theta^R_j \approx r_j+ \Theta_R, \;\; \{ i,j \} \in \{ \mathcal{B},\mathcal{R}\},  
\label{approx}
\end{equation}
where the quantities $b_i$ and $r_j$ are fluctuations with second order terms (and above) approximately equal to zero, \textit{i.e.} $b_ib_j \approx b_i r_j \approx r_i r_j \approx 0\;\; \forall \;\; i,j$. If these perturbations are sufficiently small, then $O_B \approx O_{R}\approx1$, and we have the following closed form expression for $\Delta_{BR}$ \cite{PhysicaA},
\begin{equation}
\Delta_{BR}(t) = 2\arctan\left[ \frac{C - \sqrt{\mathcal{K}} \tanh \left(\frac{t + const}{2}\sqrt{\mathcal{K}}\right) }{ \bar{\omega}_B -\bar{\omega}_R - S  } \right]
\label{AnalyticalDelta}
\end{equation}
where we have defined,
\begin{eqnarray}
\begin{split}
{\cal K} \equiv \mathcal{K}(\phi^{BR},\phi^{RB}) &=C^2 + S^2 - (\bar{\omega}_B -\bar{\omega}_R)^2, \\
C \equiv C(\phi^{BR},\phi^{RB}) &= \frac{d_T^{BR} \zeta_{BR} \cos \phi^{BR}}{|{\cal B}|} + \frac{d_T^{RB} \zeta_{RB} \cos \phi^{RB}}{ |{\cal R}|}, \\
S \equiv S(\phi^{BR},\phi^{RB}) &= \frac{d_T^{BR} \zeta_{BR} \sin \phi^{BR}}{|{\cal B}|} - \frac{d_T^{RB} \zeta_{RB} \sin \phi^{RB}}{|{\cal R}|},
\end{split}
\label{eq:ScriptS}
\end{eqnarray}
and,
\begin{equation}
d_T^{BR} = \sum_{i\in\mathcal{B}}\sum_{j\in\mathcal{R}}\mathcal{A}_{ij}^{BR}, \;\; d_T^{RB} = \sum_{i\in\mathcal{R}}\sum_{j\in\mathcal{B}}\mathcal{A}_{ij}^{RB},
\label{totaldeg}
\end{equation}
give the number of edges between Blue to Red, and Red to Blue, respectively. In the use-case of the tree
 versus the random graph $d_T^{BR} =d_T^{RB} =16$. 
Here, $const\in\mathbb{R}$ is set by the initial condition of $\Delta_{BR}$ at $t=0$.

If $\mathcal{K} \geq 0$ in Eq.(\ref{AnalyticalDelta}), \textit{i.e.} the cross couplings are sufficiently large to overcome the natural frequency difference, then $\Delta_{BR}$ approaches the stable asymptotic limit,
\begin{equation} \label{AsympDelta}
\Delta_{BR}^\infty = 2\arctan\left( \frac{C - \sqrt{\mathcal{K}}  }{ \bar{\omega}_B -\bar{\omega}_R - S } \right) .
\end{equation}
Contrastingly, if $\mathcal{K} < 0$ then $\Delta_{BR}$ remains time-dependent with periodic behaviour. 

The transient behaviour of $\Delta_{BR}$ for ${\cal K} \ge 0$ is also important for combat outcomes since, if a force can initially accelerate ahead to an early advantage even if it is unable to maintain that advantage in the asymptotic state, it may inflict enough damage to ensure victory. We use the derivative of Eq.(\ref{AnalyticalDelta}) at $t=0$ to inspect which force accelerates ahead,
\begin{eqnarray}
\begin{split}
\dot{\Delta}_{BR}(t=0) =& \frac{-\mathcal{K}}{1+\left[\frac{C-\sqrt{\mathcal{K}}\tanh\left(\frac{const \sqrt{\mathcal{K}}}{2} \right)}{\bar{\omega}_B-\bar{\omega}_R-S} \right]^2} \left[\frac{ \textrm{sech}^2 \left(\frac{const \sqrt{\mathcal{K}}}{2}\right)}{\bar{\omega}_B-\bar{\omega}_R-S} \right] \\
&= \bar{\omega}_B-\bar{\omega}_R+S,
\end{split}
\label{eq:surprising}
\end{eqnarray}
where the second line in Eq.(\ref{eq:surprising}) follows from applying the initial condition $\Delta_{BR}(t=0)=0$, leading to $const = \frac{2}{\sqrt{\mathcal{K}}}\arctanh\left(\frac{C}{\sqrt{\mathcal{K}}}\right)$. Thus, from Eq.(\ref{eq:ScriptS}) for $S$, Eq.(\ref{eq:surprising}) reveals that forces with faster mean natural frequencies or have higher cross-coupling coefficients (with the appropriate frustration) are able to gain an initial phase advantage over their adversaries leading to an early increase in their organisational combat component. Specifically, for
$\phi^{BR}=\phi^{RB}=\frac{\pi}{2}$, as in
the right-most panel of Fig.\ref{fig:CombinationPlot},
$S=0$ so that $\dot{\Delta}_{BR}(t=0)=-0.048<0$ 
giving Red a slight transient advantage.

\begin{figure}[h]
\centering
\includegraphics[width=.66\linewidth]{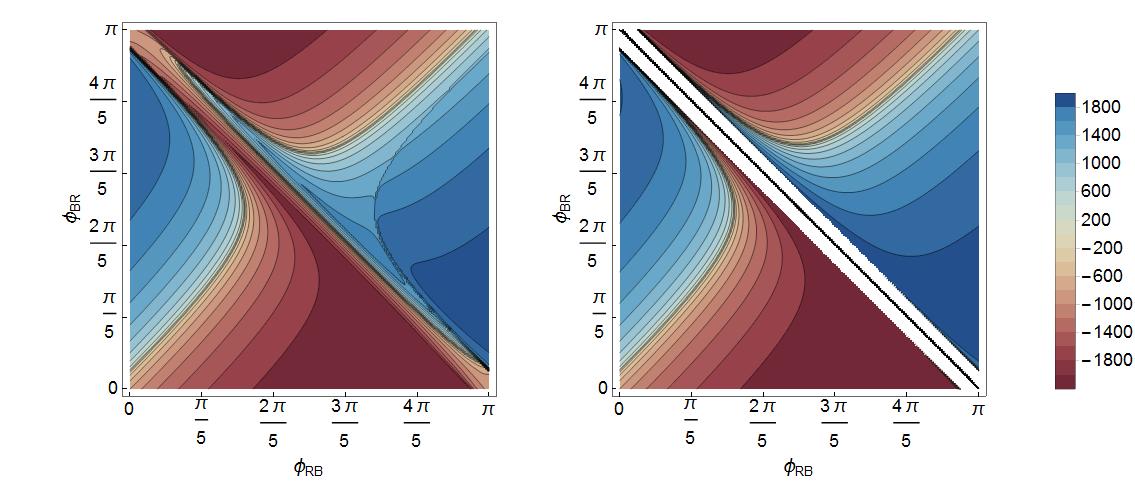}
\caption{Left most panel presents contour plot of $ p_{final}$ generated from  Eqs.(\ref{fullLan},\ref{AnalyticalDelta}). The right most panel shows the equivalent contour plot for $p_{final}$, with the asymptotic expression $\Delta^{\infty}$ which is only valid in the region ${\cal K} \ge 0$, given by Eq.(\ref{AsympDelta}).}
\label{fig:CombinationPlotanalyt}
\end{figure}

Using the same set-up as numerically
studied previously, we plot $p_{final}$ generated from the Lanchester trajectories in Eq.(\ref{fullLan}) using the assumption that $O_B \approx O_{R}\approx1$, and applying the analytic time dependent (Eq.(\ref{AnalyticalDelta})) and asymptotic (Eq.(\ref{AsympDelta})) form of $\Delta_{BR}$ in Figure \ref{fig:CombinationPlotanalyt}.
The left most panel in Fig.\ref{fig:CombinationPlotanalyt} gives the contour plot of $ p_{final}$ generated from  Eqs.(\ref{fullLan},\ref{AnalyticalDelta}), 
{\it reducing 
from 44 to 2 the dimension of the underlying system of differential equations}. Comparing with the contour plot given in Figure \ref{fig:CombinationPlot}, we see that the clustering assumption fares well in this parameter regime. 

The right most panel of Fig.\ref{fig:CombinationPlotanalyt} shows the equivalent contour plot for $p_{final}$, but with the time dependent expression for $\Delta$ in Eq.(\ref{AnalyticalDelta}) replaced with the asymptotic expression $\Delta^{\infty}$, given by Eq.(\ref{AsympDelta}), which is only valid in the region ${\cal K} \ge 0$.   We can notice that the \textit{bump} present in the left hand panel, yet absent in the right hand panel, is the effect of the transient system behaviour before the system attains steady-state. That is, if Blue locks ahead asymptotically \textit{and} accelerates ahead during the transient dynamics, it is able to obtain additional benefit.
This then explains why in the right-most case in 
Fig.\ref{fig:CombinationPlot} Red is able to win -
due to its initial transient trajectory.

\subsection{Introducing feedback into the organisational dynamics}
\label{attensec1}
\begin{figure}[h]
\centering
\includegraphics[width=1\linewidth]{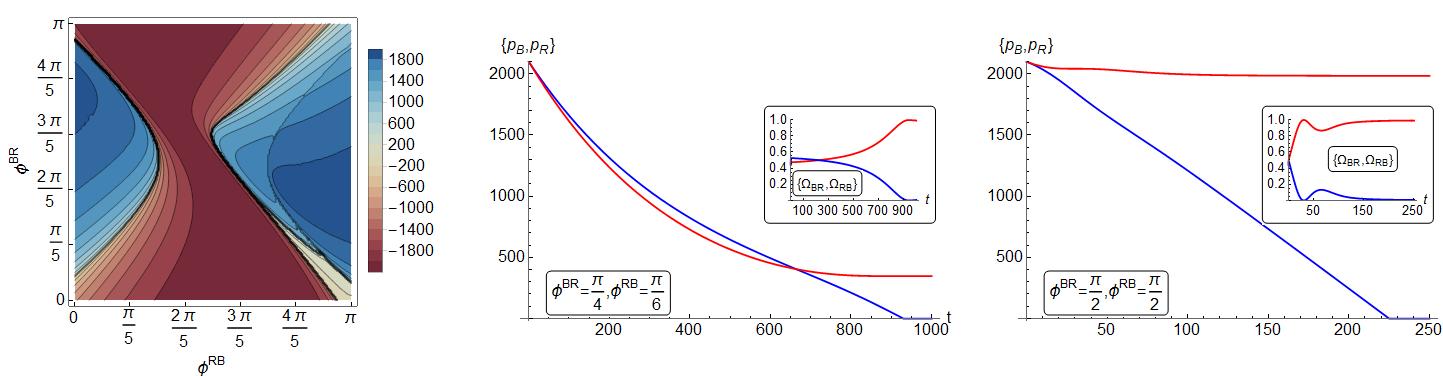}
\caption{Left most panel presents contour plot for the value of $p_{final}$ generated from  Eqs.(\ref{eq:BlueKuramoto},\ref{fullLan}), using the attenuation given in Eq.(\ref{atten1}). The middle and right-most panels give examples of the the Lanchester trajectories for specific choices of frustration parameters. Additionally, the insets of the middle and right-most plots show the dynamics with respect to time of the organisational component of each of Blue and Red's effectiveness, given in Eq.(\ref{eq:OrgComponent}).}
\label{fig:DifferentPrefactorsPlot}
\end{figure}

\begin{figure}[h]
\centering
\includegraphics[width=1\linewidth]{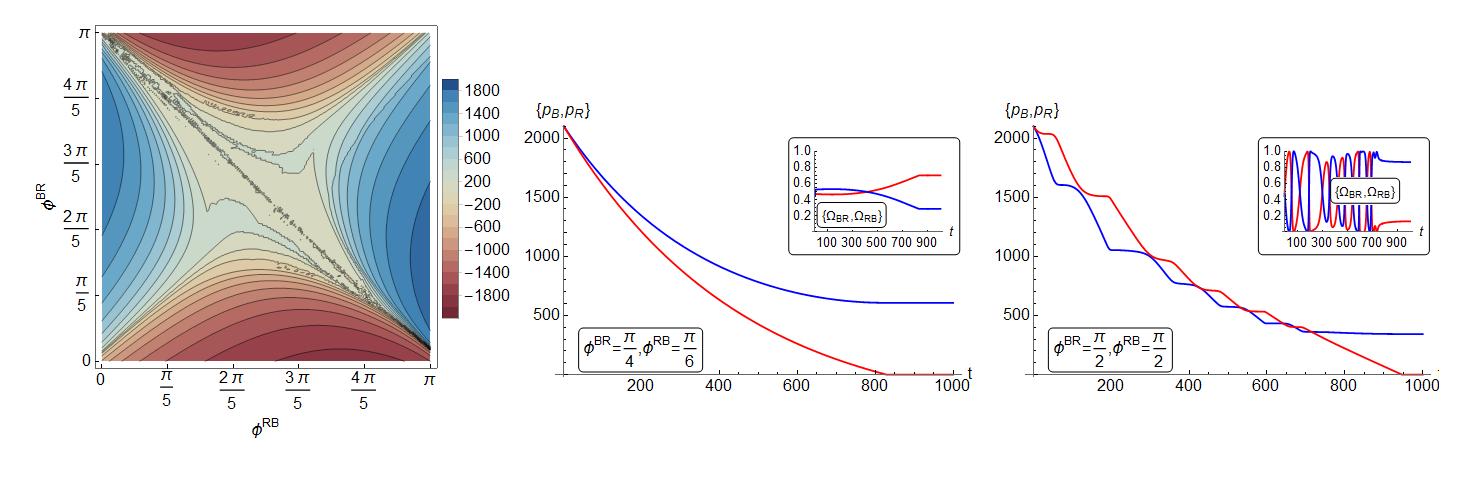}
\caption{Left most panel presents contour plot for the value of $p_{final}$ generated from  Eqs.(\ref{eq:BlueKuramoto},\ref{fullLan}), using the attenuation given in Eq.(\ref{atten2}). The middle and right most panels give examples of the the Lanchester trajectories for specific choices of frustration parameters. Additionally, the insets of the middle and right-most plots show the dynamics with respect to time of the organisational component of each of Blue and Red's effectiveness, given in Eq.(\ref{eq:OrgComponent}).}
\label{fig:DifferentPrefactorsPlot2}
\end{figure}

We now provide the most obvious generalisation of this model, where Lanchester combat outcomes feed back into the Kuramoto decision making dynamics. A straightforward way to achieve this is to attenuate the various cross couplings with a function of the current and initial populations:
\begin{eqnarray}
\begin{split}
\sigma_B &\mapsto f(p_B(t),p_{B0})\cdot \sigma_B, \quad & \sigma_R &\mapsto f(p_R(t),p_{R0})\cdot  \sigma_R,\\
\zeta_{BR} &\mapsto g(p_B(t),p_{B0},p_R(t),p_{R0})\cdot  \zeta_{BR}, \quad & \zeta_{RB} &\mapsto g(p_R(t),p_{R0},p_B(t),p_{B0})\cdot  \zeta_{RB}. \label{eq:CrsAten}
\end{split}
\end{eqnarray}
In other words, as a force suffers attrition its
ability to couple changes. 
Depending on context, when a particular force loses resource/population agents may
lose focus on staying ahead of the adversary's organisational decision cycles, in which case $g(p, p_0) = \frac{p}{p_0}$; or expend more effort on staying ahead of the adversary's decision cycles, in which case $g(p, p_0) = \frac{1}{p/p_0+\epsilon}$ for small $\epsilon>0$.
%
% The same argument may also be made to arbitrary complexity of the cross-coupling, such as, should the cross coupling depend on the relative force ratio such that:
%\begin{itemize}
%	\item smaller forces are more agile, $g(p_B,p_R) =2-\frac{2p_B}{p_R+p_B}$ ; or 
%	\item easier for well-resourced opponents to out-think, $g(p_B,p_R)= \frac{2p_B}{p_R+p_B}$.
%\end{itemize}
%In reality, these functions are likely to be complex, multifaceted, and nuanced, depending on everything from the training of forces to how their communications equipment is networked.
%
%Intuitively, these functions can be understood as either exacerbating or stabilising combat results. For example, if suffering attrition causes a force's cross coupling to decrease, it will become harder to stay ahead of the adversary's decision cycle, potentially leading to an acceleration in attrition. Conversely, if attrition drives a force to increase its cross coupling with the adversary, it becomes easier to get ahead of opposition decision making which may counteract loss due to attrition. 
Situations for different $g$ are visualised in Figs.\ref{fig:DifferentPrefactorsPlot} and \ref{fig:DifferentPrefactorsPlot2}. We also explored different values of $f$, only to produce contours very similar to Fig.\ref{fig:CombinationPlot}. We explore the reason for this in the following. 

In Fig.\ref{fig:DifferentPrefactorsPlot}, we present results for the particular choice
\begin{equation}
    f_{B} = f_R = 1, \;\; g_B = \frac{p_B}{p_{B0}}, \;\; g_R = \frac{p_R}{p_{R0}},
    \label{atten1}
\end{equation}
which represents the scenario where suffering attrition results in the degradation of the ability to stay ahead of the adversary's collective decision cycle. Focusing on the contour plot on the left hand panel, we see that the effect of the attenuation given in Eq.(\ref{atten1}) greatly expands the Red force's ability to win engagements compared to 
where attrition makes no impact on C2.
Importantly, there are now whole regions --- approximately $\phi^{RB} \in (3\pi/10 , \pi / 2)$ --- where Red wins the engagement \textit{regardless} of the choice Blue makes with its frustration parameter. This is a new situation which does not appear in the equivalent contour in Figure \ref{fig:CombinationPlot}, and highlights the difference in network topology, and the importance of a force maintaining focus on the adversary's decision cycles. 

The middle panel, with the choice $(\phi^{BR}, \phi^{RB}) = (\pi/4, \pi/6)$, presents a Lanchester trajectory where Blue is (barely) ahead for the majority of the engagement, but is overtaken at approximately $t=700$ and proceeds to extinction shortly after. The inset helps to
understand how this marked change occurs noting
the two-way influence between C2 and combat dynamics.
As attrition sets in for both forces, and therefore
degradation of coupling through Eq.~(\ref{atten1}), neither
is able to build up significant coherence
with $\Omega$ values remaining around $0.5$. Nevertheless, early in the competition ($t\approx 300$) Red
gains an albeit narrow advantage 
through a relatively
superior coherence across its members' decision states (inset). This is a consequence of its superior more
connected network structure -- again, providing
robustness. This in turn translates
into a deceleration in its attritional losses.
Even though both forces reach a point of parity
at $t=700$, the superior decision coherence
of Red at this stage means it ceases to take
losses while Blue's trajectory is unchanged. As
Blue suffers further attrition its capacity to
maintain decision coherence collapses. Red's coherence
only continues to improve because
it is able to stabilise it's locking
ahead of Blue (the attentuation is in $g$) and therefore
there ceases to be a conflict between that
and its internal drive for synchronisation.

Contrasting this with the right panel of
Fig.~\ref{fig:DifferentPrefactorsPlot}, for $(\phi^{BR}, \phi^{RB}) = (\pi/2, \pi/2)$, we see the Red force staying ahead for the entirety of the engagement. In the insets of both of these panels which show the value of the organisational factors $\Omega$ over time, we notably witness much more non-monotonicity than the equivalent insets in Fig.\ref{fig:CombinationPlot}. 

Contrastingly, in Fig.\ref{fig:DifferentPrefactorsPlot2}, we present results for the attenuation choice
\begin{equation}
    f_{B} = f_R = 1, \;\; g_B = \frac{1}{\frac{p_B}{p_{B0}}+10^{-3}}, \;\; g_R = \frac{1}{\frac{p_R}{p_{R0}}+10^{-3}},
    \label{atten2}
\end{equation}
which emulates the opposite scenario, where suffering losses results in increased effort to stay ahead of the adversary's collective decision cycle. Visual inspection of the contour plot on the left hand panel 
of Fig.\ref{fig:DifferentPrefactorsPlot2}
shows that the attenuation given in Eq.(\ref{atten2}) largely puts both forces on an even playing field;
blue and red regions of the plot are nearly equal in area.
Interestingly, we see that large regions of the $(\phi^{BR}, \phi^{RB})$ parameter space now result in stalemate scenarios for the Blue and Red forces, especially in the middle and along the diagonals. This situation is absent from Fig.\ref{fig:CombinationPlot}. From the contour in Figure \ref{fig:DifferentPrefactorsPlot2}, we see that the attenuation in Eq.(\ref{atten2}) almost eliminates the effects of
intrinsic differences between Blue and Red 
(native frequencies and/or network structure). 

The middle panel, with choice $(\phi^{BR}, \phi^{RB}) = (\pi/4, \pi/6)$, presents a fairly uninteresting Lanchester trajectory, which typifies the behaviour witnessed in the blue and red coloured regions of the contour plot with a clear winner. In real contrast to this, the right panel, for $(\phi^{BR}, \phi^{RB}) = (\pi/2, \pi/2)$, exemplifies the extent of the impact the attenuation of Eq.(\ref{atten2}) when one force gains an advantage over the other. Notably, for small times, as the difference between Red and Blue gets too large, the organisational effects (see inset) oscillate rapidly, causing the populations to achieve parity. As time evolves the behaviour repeats quasi-periodically, until Red diminishes linearly at $t\approx 700$. 

\begin{figure}
\centering
\includegraphics[width=0.66\linewidth]{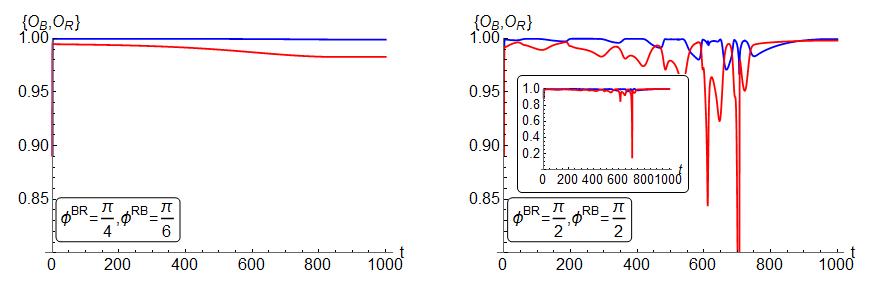}
\caption{Plots of the time dependent order parameters --- Eq.(\ref{eq:OrderParam}) --- corresponding to the Lanchester trajectories displayed in Figure \ref{fig:DifferentPrefactorsPlot2}. The inset on the right hand plot has the same parameter values as the parent panel and displays a larger range to offer more perspective.} 
\label{fig:synchanalyt}
\end{figure}

In an effort to explain this sudden drop in Red's effectiveness and lack of quasi-periodicity in the Lanchester trajectories past $t=700$ in the right-most panel of Fig.\ref{fig:DifferentPrefactorsPlot2}, we plot the corresponding time dependence of the order parameters, $O_B$ and $O_R$ from Eq.(\ref{eq:OrderParam}), for these two specific instances in Figure \ref{fig:synchanalyt}. The left hand plot, for $(\phi^{BR},\phi^{RB})=(\pi/4,\pi/6)$, displays highly synchronised Blue and Red order parameters, meaning that the multiplicative  $O$--terms contained within Eq.(\ref{eq:OrgComponent}) are approximately equal to unity, and do not negatively influence the organisational contribution to the force effectiveness. The right hand plot however shows a markedly different scenario, with both trajectories highly variable, and although the Blue order parameter stays close to unity for all $t$, the Red order parameter deviates significantly at multiple times, with the largest occurring at approximately $t=700$. This significant deviation causes the multiplicative $O_R$--term for Red's internal organisational effectiveness to plummet to zero. Moreover, we also see that Red's organisational combat factor $\Omega$ 
in the inset of the right plot of Fig.\ref{fig:DifferentPrefactorsPlot2} subsequently loses its ability to rise in a quasi-periodic manner. 
In summary, because the Red force excessively 
tightens in relation to Blue after suffering
attrition its internal decision-making coherence 
has degraded, which in turn impacts
on its combat effectiveness. Can such second-order effects be computed semi-analytically?

%Another natural extension to the model is incorporating the removal of nodes from the Kuramoto network as resource attrition occurs. Although certainly possible, it is difficult determine which nodes should be removed as the order of node removal is important to how properties such as network synchronisation are changed. For example, if a node with high betweenness centrality, such as Node 1 in the Blue Hierarchical subnetwork in Figure \ref{fig:PABlueNetwork}, then the network may be segregated into disconnected clusters that will lead to a diminishment in their global order parameter and thus combat ability. If instead each node was engaged in its own Lanchester combat with the opposing network nodes and resources were distributed and shared amongst those individual nodes, it would be simple to determine when a node should be removed. Just check if it still has any resource left! This motivates, alongside the desire to incorporate manoeuvre, a model where both combat and decision dynamics occur along a network.

\subsection{Semi-analytic approach in presence of combat feedback on organisational dynamics}
\label{semian}
\begin{figure}[h]
\centering
\includegraphics[width=0.66\linewidth]{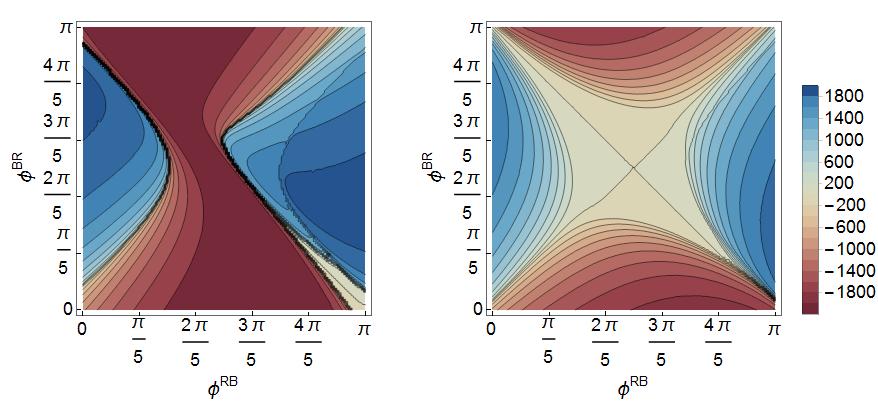}
\caption{Contour plots of the dimensionally reduced system given in Eq.(\ref{totalapprox}) with the same parameter inputs chosen to produce Figures \ref{fig:DifferentPrefactorsPlot} and \ref{fig:DifferentPrefactorsPlot2}. The left plot presents the case $g(p) = p/p_0$, which can be directly compared to the contour plot in Figure \ref{fig:DifferentPrefactorsPlot}. Additionally, the right panel applies $g(p) = (p/p_0 + 10^{-3})^{-1}$, and hence can be directly compared to the contour plot of Figure \ref{fig:DifferentPrefactorsPlot2}.}
\label{fig:analyticapp}
\end{figure}

We perform a dimensional reduction approximation to the system of Eqs.(\ref{eq:BlueKuramoto},\ref{fullLan}), in a manner equivalent to \cite{PhysicaA} in order to exploit the approximation given by Eq.(\ref{AnalyticalDelta}). 
Because Lanchester terms appear now in the
phase dynamics, the previous approach is
insufficient. 
This is a rather lengthy calculation relegated to the Appendix. The derivation exploits properties of the
graph Laplacian spectra for the Blue and Red force
networks \cite{Boll1998, PecCar1998}. Critically, 
to simplify we assume
that each force is close to 
internal synchronisation of its decision cycles,
$O_B,O_R \approx 1$.
These considerations give:
\begin{eqnarray} 
\begin{split}
\dot{\Delta}_{BR} = \bar{\omega}_B-\bar{\omega}_R -\left[ \frac{\zeta_{BR}}{|{\cal B}|} g(p_B) \sin (\Delta_{BR} - \phi_{BR}) d^{BR}_T +\frac{\zeta_{RB}}{|{\cal R}|} g(p_R) \sin (\Delta_{BR} + \phi_{RB}) d^{RB}_T \right],\\
\dot{p}_B = -\kappa_{RB} \frac{1-\sin(\Delta_{BR}) }{2}p_R\mathcal{H}(p_B), \;\;
\dot{p}_R = -\kappa_{BR} \frac{1+\sin(\Delta_{BR}) }{2}p_B\mathcal{H}(p_R).
\end{split}
\label{totalapprox}
\end{eqnarray}
In the Appendix we detail why Eq.(\ref{totalapprox}) is
independent of the feedback terms $f$, consistent
with our observations that choices for these functions
played minimal role in the contour plots. 
We note that if the remaining feedback term $g$ equals unity, then the differential equation for $\Delta_{BR}$ in Eq.(\ref{totalapprox}) is solvable, leading to the expression already given in Eq.(\ref{AnalyticalDelta}). We give contour plots of the outputs of Eq.(\ref{totalapprox}) in Fig.\ref{fig:analyticapp}, which should be compared with those in Figs.\ref{fig:DifferentPrefactorsPlot} and \ref{fig:DifferentPrefactorsPlot2}. 

From visual inspection, we see that the choice $g(p) = p/p_0$ in the left panel serves to accentuate Red's ability to win engagements, with the attenuation term lessening each force's ability to pursue their intent of staying ahead of their adversary's global phase as they suffer attrition. Contrast to this to the right most-plot in Fig. \ref{fig:analyticapp}, we see that the choice $g(p) = (p/p_0 + 10^{-3})^{-1}$ largely levels the ability for Red and Blue to win engagements.  However, as found in Fig. \ref{fig:synchanalyt}, 
the tightening of a force's attention on the adversary
that this models can become out of balance for the 
requirement of internal coherence, as witnessed in Sec. 4.4 of \cite{PhysicaA} --- explaining the difference in detail between the right panel in Fig.\ref{fig:analyticapp}, and the equivalent contour plot in Fig.\ref{fig:DifferentPrefactorsPlot2}. Notably, in Fig.\ref{fig:analyticapp} we are not witnessing any internal network decoherence of the force as the approximations 
used in the Appendix for Eq.(\ref{totalapprox})
exclude this behaviour.

To summarise this section: using
a representation of homogeneous forces,
though networked in their decision-making, we have established the principles for unifying the Lanchester and Kuramoto
models. We did this for both cases where
attrition does not, or does, impact on the abilility of the forces
to conduct C2. Additionally, we see that, with some
considered approximations, we are able to reproduce, and
analytically understand, many 
of the predictions of the models for the outcome of
battle. We now take these principles and insights
to highly non-homogeneous forces.

\section{Embedding Lanchester into a network: the networked Lanchester-Kuramoto-Sakaguchi model}
\label{sec3}
We finally present the model we have been
building up to, a representation for two distributed inhomogeneous forces in combat, where their internal
decision-making dynamics is in tension with
their seeking to achieve tactical decision-making ahead of the
adversary. To the degree that they have decision advantage
they gain a combat advantage. Moreover, they
have the capacity to shift force elements dynamically
through a network adapting to the dynamics of battle,
where we build in the Lanchester-manoeuvre parts of \cite{Kall19}.
%We define the model, consider it the same use-case as applied to the global model, and then to a very stylised inhomogeneous network. This elucidates when manoeuvre dominates over attrition. We then move to a less symmetric, pseudo-realistic scenario. We show the extent to which aspects can be analytically understood using techniques developed earlier.

\subsection{Model definition}
\label{sec3.1}
To generalise the model given in Sec.\ref{sec2} by placing the Lanchester component on a network, introducing manoeuvre, we propose the following 2$|\mathcal{K}|$--dimensional system,
\begin{eqnarray}
\dot{\theta}_i =\mathcal{H}(p_i)\left[\omega_i-\sum_{j\in\mathcal{K}} \mathcal{H}(p_j)\mathcal{K}_{ij}\sigma_{ij}\sin(\theta_i-\theta_j-\Phi_{ij})\right], \;\; i \in \mathcal{K}, \label{NetKura}\\
\dot{p}_i = \mathcal{H}(p_i)\left[ \sum_{h\in\mathcal{M}}\mathcal{H}(p_h)\mathcal{M}_{ih}\frac{\Gamma_{i,h}+\Gamma_{h,i}}{2}(\delta_h p_h - \delta_i p_i)\frac{\cos(\theta_h-\theta_i)+1}{2} \right. \nonumber\\
-\left. \sum_{k\in\mathcal{E}}\mathcal{H}(p_k)\mathcal{E}_{ik}\kappa_{ik}p_k\text{d}_k \frac{\sin(\theta_k-\theta_i)+1}{2}O_k, 
\right],\;\; i \in \mathcal{K}. \label{NetLan}
\end{eqnarray}
The organisational component Eq.(\ref{NetKura}), is largely unchanged from Eqs.(\ref{eq:Kuramoto},\ref{Kuradef}) apart from additional Heaviside factors that ensure 
force elements reaching zero (or some other threshold) are
incapable of further participating in the dynamics.
The need for this arises here
specifically because of the localised
interaction between decision states (Kuramoto phases) $\theta_i$
and individualised combat element force strengths
$p_i$. Quite simply, if an individual force element is
extinguished it ceases to couple to
partners -- the {\it internal} factor ${\cal H}(p_i)$ in
Eq.~(\ref{NetKura}) -- {\it and} it ceases to
cycle through its own disconnected decisions
-- the {\it external} factor ${\cal H}(p_i)$ in
Eq.~(\ref{NetKura}).

Focusing on the networked Lanchester terms Eq.(\ref{NetLan}), we distinguish between \textit{manoeuvre} $\mathcal{M}$ and \textit{engagement} $\mathcal{E}$ components of the network $\mathcal{K}$:
\begin{equation}
\mathcal{K} = 
\underbrace{\left[ 
\begin{array}{cc}
0 & \mathcal{A}^{BR} \\ 
\mathcal{A}^{RB} & 0
\end{array} \right]}_{\equiv \mathcal{E}}
+
\underbrace{\left[ 
\begin{array}{cc}
\mathcal{B} & 0 \\ 
0 & \mathcal{R}  
\end{array} \right]}_{\equiv \mathcal{M}}.
\label{EandM}
\end{equation}
Additionally, we define the \textit{local order parameter} $O_k$, at each node $k$, via \begin{equation}
   O_k =   \frac{ \left|\sum_{m\in{\cal{M}}}{\cal M}_{km}\mathcal{H}(p_m)e^{\sqrt{-1} \theta_m} \right| +\epsilon_2}{\sum_{m\in{\cal{M}}} {\cal M}_{km}\mathcal{H}(p_m) + \epsilon_2}  , 
   \label{localord}
\end{equation}
which measures how synchronised each oscillator is with its nearest neighbour. In contrast to $O_B, O_R$ used
previously, the normalisation here is time-dependent, and may vanish, hence the regularisation
$\epsilon_2$. It should be clear how this insertion
of $O_k$ in the last term of Eq.(\ref{NetLan}) with the sine factor, generalises the principles in Eq.(\ref{fullLan}).

The parameters $\delta$ and d in Eq.(\ref{NetLan}) are defined by,
\begin{equation}
\delta_k =\frac{1}{\sum_{m\in\mathcal{E}}\mathcal{E}_{km}p_m+ 1}, \;\; \text{d}_k = \frac{1}{\sum_{m\in\mathcal{E}}\mathcal{E}_{km}\mathcal{H}(p_m)+\epsilon_2},
\label{deltad}
\end{equation}
introduced in \cite{Kall19}, with $\delta$ moderating the flow of resource in the internal networks, and $\text{d}$ moderating the effectiveness of the engagement. We reiterate that $\epsilon_2$ is a small real number, chosen in this work specifically by $\epsilon_2=10^{-20}$.

%The first line of Eq.(\ref{NetLan}) is the main deviation from  Eq.(\ref{fullLan}). 
Eq.(\ref{NetLan}) draws upon the
model developed (without C2) in \cite{Kall19}. The differences of this extension are characterised by
distributing resource across network nodes rather than a global total; dividing the combat effectiveness of a node by the number of concurrent engagement connections $\text{d}_k$;
replacing the network-average phase difference term $\Delta_{BR}$, with a nodal nearest neighbour adversary factor, $\theta_k-\theta_i$; and
replacing the global order parameter $O_{B}$ and $O_{R}$ by a local per-node check against any other nodes connected along the manoeuvre network, $O_{k}$, defined in Eq.(\ref{localord}).

The manoeuvre mechanism given in the first line of Eq.(\ref{NetLan}) allows node $h$ with non-zero resource and manoeuvre connection to node $i$ ($\mathcal{M}_{ih}=1$) to share resource
according to the average of their manoeuvre constants $(\frac{\Gamma_{ih}+\Gamma_{hi}}{2})$,
multiplied by the difference in force ratios between the sharing nodes and the total adversarial population among their respective engagements $(\delta_h p_h - \delta_i p_i)$,
and the extent to which their decision making is in agreement $(\frac{\cos(\theta_h-\theta_i)+1}{2})$. Alternately,
$\mathcal{M}$ may be thought of as a network regulating
the force-flow of through the combat engagement.

The choice of the term `1' in the denominator of $\delta$ in Eq.(\ref{deltad}) represents a `standing force' that a non-combat node wishes to maintain. If replaced by a number less than unity, non-combat nodes give up their resource rapidly to engaged nodes; if greater than unity, all nodes distribute resource evenly with little regard for engagements.

The matrix terms $\kappa$ and $\Gamma$ control the rates of physical effectiveness and manoeuvre, respectively:
\begin{equation}
\kappa = \left[ 
\begin{array}{cc}
0 & \kappa_{BR} \\ 
\kappa_{RB} & 0
\end{array} \right], \;\; 
 \Gamma = \left[ 
\begin{array}{cc}
\gamma_B & 0 \\ 
0 & \gamma_R
\end{array} \right].
\end{equation}

It is important to note that the choice of manoeuvre using the difference in force ratios is a heuristic. Indeed, many historical engagements have been decided by the ability to rapidly concentrate and distribute forces as required. 
In \cite{Kall19}, it was seen that optimised networks (without an explicit
dynamical representation of C2) with such a 
heuristic
delivered structures that resembled Manoeuvre Theory
concepts, such as Feints and Disruptive Fire.

The current choice for $\delta$ in this work disperses resources among engagements to match the opposition, rather than exploiting weaknesses with targeted force concentrations due to application of human initiative. 
Nevertheless, one can argue that forces behave in a somewhat diffusive way in various contexts. Although manoeuvres resulting in localised concentrations often require knowledge of the state of combat beyond the nearest neighbour,  if agents can only act on local information 
{\it in lieu} of a network-wide strategy, this may be the most optimal choice available. Moreover, particularly risk-averse decision-makers may still order this as a global strategy despite being sub-optimal in many scenarios as seen by the French Military ``penny packet'' distribution of superior armour in the prelude to the WWII Battle for France. At any rate, we save the more advanced heuristics, such as adversarial dynamic reallocation of forces, to future studies.

A final difference between the global  (Eqs.(\ref{eq:Kuramoto},\ref{fullLan})) and networked (Eqs.(\ref{NetKura},\ref{NetLan})) Lanchester models 
is that when 
nodes are extinguished
we set the corresponding {\it phase velocity} to zero; `death' of the combat element freezes
its decision-making capability.
While not strictly necessary, as empty nodes cannot affect their nearest neighbours, it does have the advantage of reducing the dimensionality of the dynamical system as time progresses which slightly improves computation time.

\subsection{Pairwise Lanchester feedback and numerical comparison}
\label{secatten2}
Similar to Sec.\ref{attensec1}, we now propose the following attenuation functions for the feedback of Lanchester dynamics into the Kuramoto-Sakaguchi component of the model in Eq.(\ref{NetKura})
\begin{eqnarray}
\begin{split}
\sigma_B &\mapsto \frac{2p^B_j}{p^B_i+p^B_j}\cdot \sigma_B, \quad & \sigma_R &\mapsto \frac{2p^R_j}{p^R_i+p^R_j}\cdot  \sigma_R,\\
\zeta_{BR} &\mapsto \frac{2p^B_i}{p^B_i+p^R_j}\cdot  \zeta_{BR}, \quad & \zeta_{RB} &\mapsto \frac{2p^R_i}{p^R_i+p^B_j}\cdot  \zeta_{RB} 
\end{split}\label{eq:crossFeedback}.
\end{eqnarray}
The attenuations multiplying the intra-network couplings $\sigma$ here mean that nodes with higher resource have more weight over their nearest neighbours' (in the manoeuvre network) decision making, and are able to track an engaged adversary better if those respective neighbours and adversaries have less resource. The factor of `2' in the numerator ensures that when resources are equal, the resource coefficient is unity. Also, the attenuations multiplying the couplings $\zeta$, as with the equivalent in Eq.(\ref{atten1}), sharpen the transitions which reinforce early winning trajectories. Before detailing the use-cases of the networked model defined in Sections \ref{sec3.1}--\ref{secatten2}, Table \ref{tab:tab2} presents the parameter values applied for each of the use-cases.

\begin{table}[ht]
\caption{Summary of parameter values used for the three use-cases in Sections \ref{use1}--\ref{use3}.} % title of Table
\centering % used for centering table
\begin{tabular}{c c c c} % centered columns (4 columns)
\hline %inserts horizontal line
 parameters & use-case 1 & use-case 2 & use-case 3 \\ [0.5ex] % inserts table
%heading
\hline % inserts single horizontal line
$\{\sigma_B,\sigma_R\}$ & $\{8,0.5\}$ & $ [0,5]$ & $\{[0,0.5],0 \}$\\
$\{\bar{\omega}_B, \bar{\omega}_R\}$ & $\{0.503,0.551\}$ & 1 and 1.009 & \{1.104,0\} \\
$\{\zeta_{BR}, \zeta_{RB}\}$ & 0.4 & $\{0,0.1,0.5,1,1.5,2\}$ & 0.5 \\ && and $\{0,1.25,3\}$ \\
$\{\kappa_{BR},\kappa_{RB} \}$ & 0.005 & 0.1 & $\{ [0,0.03],0.01\}$ \\
$\{\phi^{BR},\phi^{RB}\} $ & $ [0,\pi]$ & $\frac{\pi}{4}$ & $\{\frac{\pi}{4},0\}$ \\ 
$\{\gamma_B, \gamma_R\}$ & 1 & 1 & $\{\{1,5,10,15\},0\}$\\
[1ex] % [1ex] adds vertical space
\hline %inserts single line
\end{tabular}
\label{tab:tab2} % is used to refer this table in the text
\end{table}

\subsection{Numerical calculations with use-case 1}
\label{use1}
Due to its complexity, this model must be solved numerically. However, to enable potential understanding of underlying mechanisms for model behaviour, we additionally look at the model through the lens of the semi-analytic approximation, detailed in Section \ref{semian}, when appropriate. 
To facilitate the embedding of manoeuvre here, the $\mathcal{E}$ and $\mathcal{M}$ adjacency matrices in Eq.(\ref{EandM}) were respectively introduced to define engagement and manoeuvre network connections. Although these are fundamentally different, we simplify by identifying Kuramoto intra-network connection with manoeuvre connections, and Kuramoto inter-network connection with engagement/combat connections.
We thus initially explore the same use-case for the global model.

\begin{figure}[h]
\centering
\includegraphics[width=0.7\linewidth]{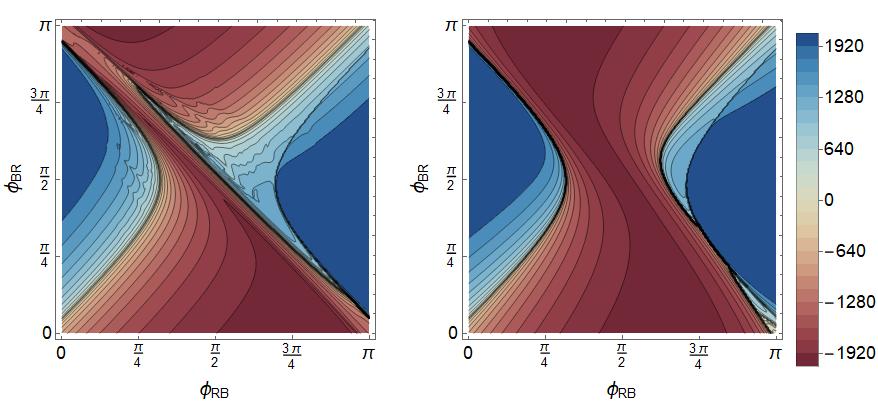}
\caption{Contour plots of combat outcomes against frustrations $\phi_{RB},\phi_{BR}\in \mathbb{S}^1$ for use-case 1. Left and right panels give Network Lanchester (Eq.(\ref{NetLan})) outcomes without and with cross-coupling feedback of the form in Eq.(\ref{eq:crossFeedback}), respectively.}
\label{fig:PhysicaAPlots}
\end{figure}

In Fig.\ref{fig:PhysicaAPlots} we give plots of Eq.(\ref{NetLan}), which should be compared to the contours in Figs.\ref{fig:CombinationPlot} and \ref{fig:DifferentPrefactorsPlot} for the global model. We run the networked model on the same scenario and parameter values as described in 
Sec.\ref{sec2} to provide comparison. With
the generalised model the additional parameters are set as follows: the new manoeuvre constant is set to $\Gamma=1$ between all manoeuvre connections; the total population of 21
nodes for each network is evenly distributed 100 resources per node for both Blue and Red. Lastly, we generalise the expression of $p_{final}$ via,
\begin{equation}
p_{final} =  \sum_{i \in \mathcal{B}} \left. p_i^B(t) \right|_{t=t_{final}} -  \sum_{i \in \mathcal{R}}  \left. p_i^R(t) \right|_{t=t_{final}}.
\label{pfinal2}
\end{equation}
from Eq.(\ref{pfinal1}) to account for the fact the resource is now distributed through each network.

Overall, visual inspection reveals that both panels in Fig.\ref{fig:PhysicaAPlots} which display the outputs of the networked Lanchester model are very similar to their global counterparts in Figs.\ref{fig:CombinationPlot} and \ref{fig:DifferentPrefactorsPlot}. For scenarios which exhibit a relatively high degree of intra-network synchronisation, as well as a homogeneous distribution of resource and engagement, the networked Lanchester model in Eq.(\ref{NetLan}) essentially reduces to the global formulation given in Eq.(\ref{fullLan}). This is to be expected given that homogeneity removes the necessity to manoeuvre resource as all force ratios are equal from the start of combat. Attrition dominates.
This result justifies the application of the global Lanchester analysis, using Eq.(\ref{fullLan}), to understand behaviours seen in the networked Lanchester expressions of Eq.(\ref{NetLan}), when in the appropriate regimes. However, as we are interested in exactly those scenarios where attrition doesn't dominate over
manoeuvre, we now consider a stylised inhomogeneous case.

\subsection{Use-case 2: The `Fighting Fish'}
\label{fightfish}
\begin{figure}[h]
\centering
\includegraphics[width=0.3\linewidth]{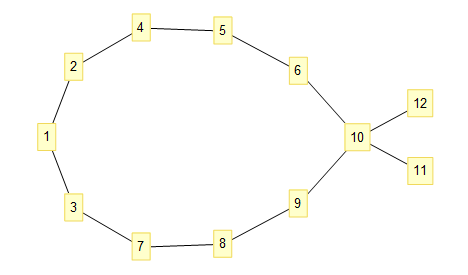}
\caption{Graph diagram of the Blue and Red networks, $\mathcal{B}$ and $\mathcal{R}$.}
\label{fig:FFNet}
\end{figure}

To highlight the nuances of the networked model, we consider a heterogeneous scenario where the trade-off between manoeuvre and engagement effectiveness has no analogue in the global Lanchester model. We construct a network and population distribution in which the ability to manoeuvre is implicitly rewarded. Both the Blue and Red networks, $\mathcal{B}$ and $\mathcal{R}$, consist of a ring with a hub node that is connected to two additional nodes as depicted in Figure \ref{fig:FFNet}. Combat occurs in one-to-one engagements across the first three nodes, \textit{i.e.} the engagement adjacency matrix $\mathcal{E}$ consists of Blue node $j$ connected to Red node $j$, for $j \in \{1,2,3\}$.
Viewing these two graphs with mutually combat nodes adjacent to each other suggests the name
`fighting fish'.

\begin{figure}[h]
\centering
\includegraphics[width=0.93\linewidth]{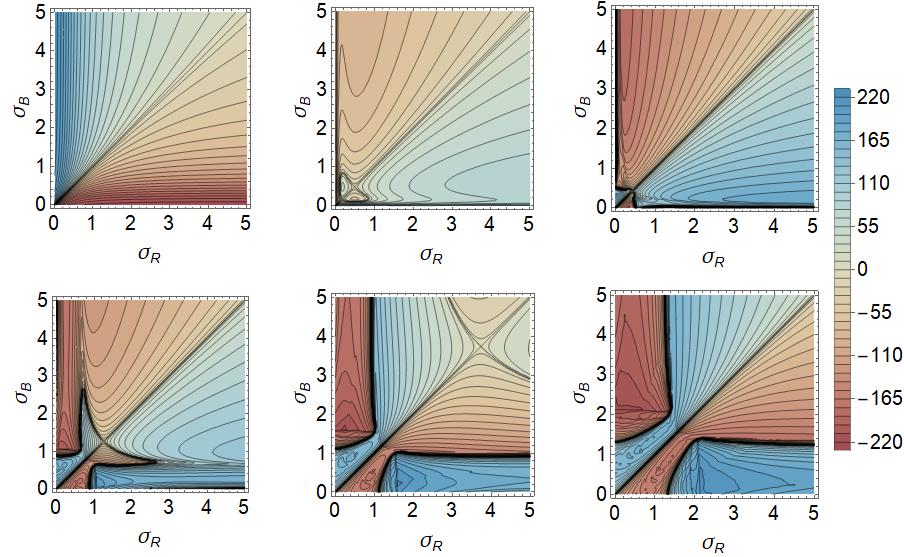}
\caption{Engagement outcomes for Eq.(\ref{NetLan}) with networks networks $\mathcal{B}$ and $\mathcal{R}$ given use-case 2 in Fig.\ref{fig:FFNet}. Plots detail the effect of intra-network couplings over the domain $(\sigma_R,\sigma_B)\in[0,5]\text{x}[0,5]$, with a simultaneous increase of inter-network couplings $\zeta_{RB} = \zeta_{BR} \equiv \zeta$ from left to right given by $\zeta= \{0,0.1, 0.5\}$ for the top row, and $\zeta= \{1,1.5, 2\}$ for the bottom row.}
\label{fig:FFCPP}
\end{figure}

Initially, resources are distributed unevenly, with $p_1=p_2=...=p_{10} = 10$, and $p_{11}=p_{12}=100$, most on the opposite side of the engagement point. The force capable of shifting these reserves to the engagement faster than the opposition fights achieves advantage. We set the frustrations $\phi_{RB}=\phi_{BR}=\frac{\pi}{4}$, the manoeuvre constants $\gamma_B=\gamma_R=1$, and the lethality constants $\kappa_{RB} = \kappa_{BR} = 0.1$, and vary the intra-network internal couplings $\sigma_B$ and $\sigma_R$, with both inter-network couplings $\zeta_{BR}=\zeta_{RB} = \zeta$. We equate the natural frequencies, $\omega_B=\omega_R=1$, and equi-space the initial node phases within $[-\frac{\pi}{4},\frac{\pi}{4}]$; we include the feedback Eq.(\ref{eq:crossFeedback}) so capability to rapidly manoeuvre provides a benefit to adversarial decision making. Thus any resulting combat advantage follows from coupling and network structure. 
%The initial phases allow internal coupling potential to improve the initial incoherence. 
We evolve the system to $t_{final}=10^3$ and plot Eq.(\ref{pfinal2}) for various $\zeta$ in Fig.\ref{fig:FFCPP}.

Despite the structural
symmetry between the forces we see in Fig.\ref{fig:FFCPP} 
a rich picture as inter-network coupling increases, with multiple regions and transitions from competing engagement and organisational effects. 
Across the panels of Fig.\ref{fig:FFCPP} we observe the propagation of particular front.

As inter-network coupling $\zeta_{BR}$ and $\zeta_{RB}$ increases, the phase oscillators attached to engagement nodes synchronise less and less with their respective network's non-combat nodes, while potentially frequency synchronising with their corresponding adversary. Conversely, increasing intra-network coupling $\sigma_B,\sigma_R$ 
has the tendency to enable all nodes within Blue or Red to synchronise collectively.
These considerations explain the diagonal symmetry of the plots.

In detail, for $\zeta = 0$, nodes within each network experience no tension from inter-network interaction. Increasing intra-network coupling here gives for each a larger local order parameter and nearest-neighbour agreement, boosting combat effectiveness and speed of resource reallocation.
At $\zeta = 0.1$, with inter-network coupling, engagement nodes start to achieve a phase advantage over their respective adversaries. If Red has zero intra-network coupling and Blue slightly increases its own ($\sigma_R = 0, \sigma_B=0.5$), Red engagement nodes achieve a phase advantage over Blue engagement nodes,
which are now excessively seeking
internal coherence over more rapid manoeuvre. If Red increases intra-network coupling ($\sigma_R = 0.2, \sigma_B=0.5$), it loses phase advantage from its own 
excessive effort at internal synchronisation. 
These effects increase with larger $\zeta$.
Overall, the basic reason for the changing
structures within a diagonal region is the relative imbalances between decision-speed, resource-manoeuvre speed, and combat speed
resulting in transitions from Blue
combat advantage to disadvantage.
In Appendix B we examine a path through
the landscape to understand more deeply
the different variations in these imbalances.

\begin{figure}[h]
\centering
\includegraphics[width=0.93\linewidth]{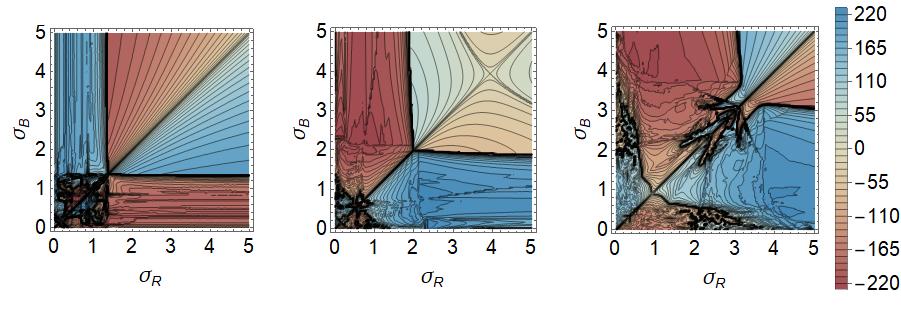}
	\caption{A selection of phase plots with increasing cross coupling $\zeta = \{ 0,1.25,3\}$ from left to right
	for use-case 2 with natural frequencies randomly chosen from a uniform distribution over $[0,2]$.}
	\label{fig:FFCPPRand}
\end{figure}

But thus far, with symmetric networks
and parameters we
have obtained quite regular structures
in the landscape. To see a first effect of randomness on this, we sample natural frequencies from a uniform distribution between $[0,2]$, with a selection of resulting phase plots shown in Fig.\ref{fig:FFCPPRand}.
We see here remnants of the order
in Fig.\ref{fig:FFCPP} but with
elements of disorder. The left panel in Fig.\ref{fig:FFCPPRand} shows that even for zero inter-network coupling the regions of differentiated behaviour are clearly present as now internal coupling must overcome the tension created by the heterogeneous natural frequency distribution. 
%This also manifests itself in the plots with non-zero $\zeta=1.25$ in the middle panel, as the purely dynamic phased oscillator region is a larger subset of the domain than in the uniform natural frequency case for Figure \ref{fig:FFLabelledPhasePlot}; a symptom of the larger stressors/requirements for each network to achieve internal synchronisation. We note also that similar structure to the uniform natural frequency cases in Figures \ref{fig:FFCPP} and \ref{fig:FFLabelledPhasePlot} exist in the top right 'steady' regions. 
However, the regions in Fig.\ref{fig:FFCPPRand} along diagonals, where there is less
likelihood of one force or the other achieving 
internal synchronisation, display considerable complexity. This suggests that for parameter
settings excluding 
steady-state the outcome is more sensitive to heterogeneous natural frequency distributions.
%
%We can also observe the propagation of this region of complexity by inspecting the combat results as a function of frustration. We set $\sigma_B=\sigma_R=1.5$ and sweep through cross couplings plotting Frustration of Blue on Red $\phi_{BR}$ and Red on Blue $\phi_{RB}$ against end of combat population.
%
%\begin{figure}[h]
%	\centering
%	\includegraphics[width=0.8\linewidth]{FFFrusCrosstot.png}
%	\caption{A sequence of plots investigating the interaction of frustrations $\phi_{RB},\phi_{BR}\in[0,\pi]$ with a simultaneous increase of both cross couplings $\zeta_{RB},\zeta_{BR}$ from top left to bottom right.}
%	\label{fig:FFCPP}
%\end{figure}
%
%
%We see that and internal coupling of $\sigma = 1.5$ is sufficient to maintain order initially, however, complex structures develop once the cross coupling is of comparable magnitude. Note that by $\zeta=2$, the point where $\sigma_B=\sigma_r=\frac{\pi}{4}$ is just on the edge of complexity in the landscape as we would expect from studying Figure \ref{fig:FFCPP_C2}. Given the frustration setting directly affects the strength of mutual engagement frequency forcing, we are observing how more ambitious settings for frustration, particularly when both subnetworks are greedy, can trigger instabilities at lower cross coupling values. 
%
Next we consider a more heterogeneous network structure.

\subsection{Use-case 3: complex network scenario for Blue.}
\label{use3}
Consider now a Blue force spread across a 50 node network. Here, each of Blue's nodes are resource hubs, with the connections representing routes through which
resource may be transported. Contrastingly, the Red force is not networked with 5 detached nodes
each of which engage a single Blue node. Thus,
linking the equivalently indexed $\{1,\dots,5\}$ nodes of
Blue and Red gives the engagement network $\mathcal{E}$. 
We pit Blue's extreme advantage of being networked against a Red overwhelming
force concentration. Both Blue and Red possess a total resource $13750$. Red evenly divides its force across its detached elements at 2750 per node.
Blue has 500 per node of its engagement network, and 250 in each of the remaining nodes.
Blue must thus manoeuvre reserves to the first five nodes while ensuring that they are not exhausted immediately by Red's initial superior numerical advantage. Can superior decision
making and manoeuvre overcome
numerical disadvantage? 
\begin{figure}[h]
\centering
\subfloat[][Blue network]{\includegraphics[width=0.5\linewidth]{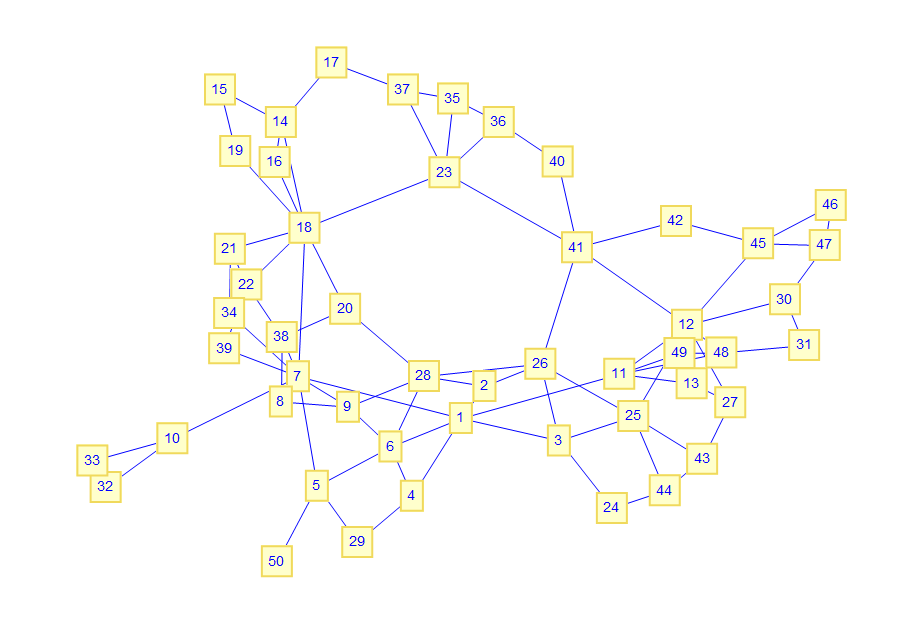}}
~
\subfloat[][Red network]{\includegraphics[width=0.2\linewidth]{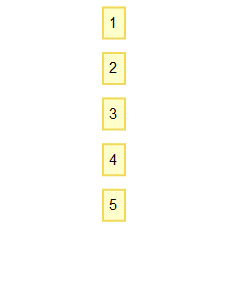}}

\caption{A graphical depiction of the network structure for each force. Blue is generated from $N_B= 50$ nodes, connected via transportation routes, and the Red network is a series of $N_R=5$ disconnected nodes. Nodes labelled $\{1,\dots,5\}$ of each network are connected with their equivalently numbered adversarial counterparts, and hence form the engagement network $\mathcal{E}$. }
\label{fig:NetworksPlotAlien}
\end{figure}

Setting $\phi_{BR} = \frac{\pi}{4},\ \zeta_{BR}=\zeta_{RB}=0.5,\ \kappa_{RB} = 0.01$, we also set $\phi_{RB}=0$ to simulate an adversary unaware of the benefits of decision advantage. We 
choose an instance of natural frequencies for both Blue and Red from a uniform distribution in the range $[0,2]$. With the Red force a set of disconnected nodes, we may set $\sigma_R = 0,\  \gamma_R = 0$. 

\begin{figure}[h]
\centering
\includegraphics[width=0.7\linewidth]{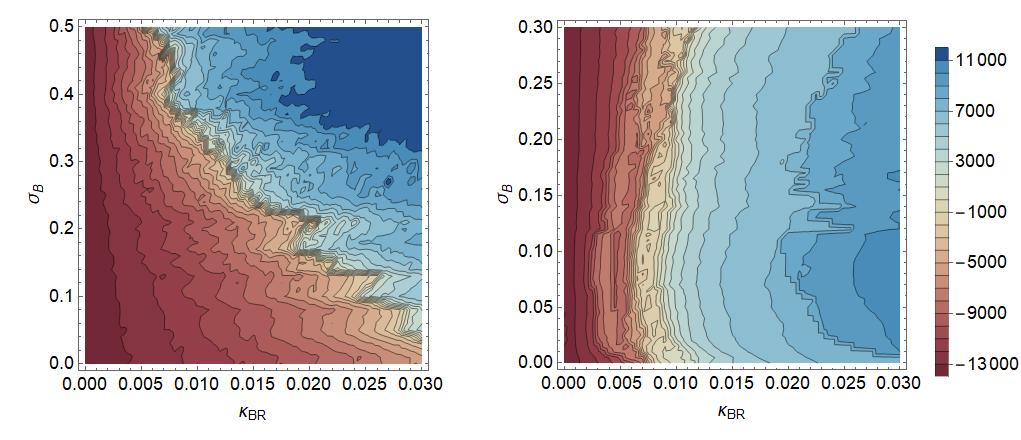}
\caption{Comparison of the global Lanchester model of Eq.(\ref{fullLan}), with the networked counterpart --- $\gamma_B=1$ --- for use-case 3.}
\label{fig:finalResults1}
\end{figure}

In Fig.\ref{fig:finalResults1} we compare the global model with manoeuvre rate $\gamma_B =1$ with the result of sweeping through $\kappa_{BR} \in [0, 0.03]$ and $\sigma_B \in [0, 0.5]$ in the networked model.
Visual inspection 
%of Fig.\ref{fig:finalResults1} 
shows that increasing internal coupling for Blue is significantly more beneficial in the global than in the networked model. This is due to the global order parameter factor in the former: it penalises systems that are unable to achieve whole of network synchronisation, namely the disconnected Red force, giving Blue a large advantage at high internal couplings. 
In the networked case, each disconnected Red node automatically has a local order parameter of unity
so now Red 
does not suffer from internal decision 
incoherence (at low coupling) attenuating its lethality; correspondingly
Blue at too
low a coupling is not drawing upon the advantages of better connectivity.
Compared to the idealised use-case 2,
Fig.\ref{fig:finalResults1} shows no symmetry
around a diagonal because the networks are so different. 
With sufficient internal coupling, Blue is able 
to overcome Red's numerical strength by swift
decision-making to manoeuvre its resources
into the engagement.

\begin{figure}[h]
\centering
\includegraphics[width=.95\linewidth]{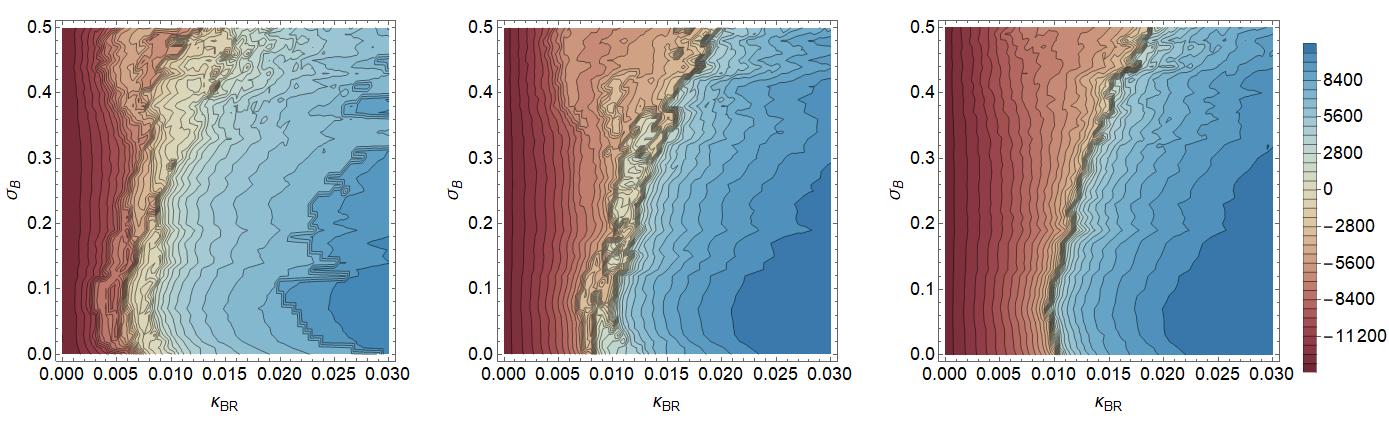}
\caption{Plot of networked Lanchester models on use-case 3, with the manoeuvre rate increasing $\gamma_B = \{ 5,10,15\}$ from left to right panels.}
\label{fig:finalResults2}
\end{figure}

Varying the manoeuvre rates to
$\gamma_{B} = \{5,10,15\}$, gives Fig.\ref{fig:finalResults2}.
The structures seen in each panel %of Fig.\ref{fig:finalResults2} 
resemble earlier analysis but with more irregularity due to
the heterogeneity in use-case 3.
As in use-case 2, we see that a degree of intra-network coupling  enables the transfer of resource to engagement nodes. 
{\it Too much coupling
generates diminishing returns}, seen in the blue-red boundary
eventually curving away
from the vertical axis.
High
internal coupling reflects an
introspective Blue, unresponsive to adversarial decision cycles, with inferior combat power. 

%This is reminiscent of the \textit{Resource Allocation Syndrome} from management literature \cite{Engwall03, Zika06}, where organisations become overwhelmed by internal resource prioritisation (personnel, funding \textit{etc}.); the organisation finds itself unable to respond to the external environment.

Lastly, an increase in the Blue manoeuvre parameter is not strictly advantageous for Blue's outcome landscape:
across the values of $\gamma_B$ there is no
sharp increase of blue regions, only
smoothening of the boundary.
This suggests a subtle interplay between manoeuvre and network structures (which
have not been optimised,
unlike in \cite{Kall19}) that warrants further exploration.

%\subsection{Interpretation in an
%Operational Research study}
Albeit with a fictitious use-case, Fig.\ref{fig:finalResults2} allows
examination of trade-offs in investment
of three types of military systems:
weapons, mobile transport, and communications infrastructure
which have proxies in the 
lethality, manoeuvre and internal coupling constants.
Against a specific scenario for the adversary
and their weapons,
and the physical layout of where Blue might intend to confront them
we observe a small `sweet spot', where
the Blue force trade-offs across these three
expenses provide advantage, at slightly lesser
lethality than the Red force weapons systems. Naturally, the recommended
trade-off point between these investment costs
needs also to be tested for sensitivity to variations in choices of
natural frequencies
(decision-making speeds) and
initial conditions, quite straightforward in models of this form.
Certainly, a force able
to {\it arbitrarily} increase the lethality of its
weapons eventually always wins over
an agile manoeuvring opponent: 
further horizontal shifts to
the right by Blue in Fig.\ref{fig:finalResults2}
eventually cross
the boundary for victory.
This captures the adage,
attributed to Stalin, that
quantity is its own quality.
%However, it may be more expensive to move horizontally through the landscape than vertically or across the panels of Fig.\ref{fig:finalResults2}. Thus {\it costs} may change the calculus and significance of a `sweet spot'.

\section{Conclusion}
We have proposed a mathematical model that unifies
the warfighting functions of
C2 and Manoeuvre
with that of Fires --- using the celebrated
Lanchester equations as substratum.
The model draws upon
multi-layered networks with different structures
available for distributed decision-making,
resource manoeuvre, and 
enemy engagement, thus exploiting aspects of modern complex systems theory.
Though deterministic,
the model exhibits adaptive dynamics with inhomogeneous forces. 
For more homogeneous scenarios,
semi-analytical approaches may be taken.
With inhomogeneity, the model remains parameter-parsimonious, allowing for efficient computational data farming.

The opportunities in future work are considerable. Optimisation of the networks, similar to \cite{Kall19}, should provide insights
into `tightness' and centralisation of decision-making to maintain agility of combat forces. Different heuristics for how force elements transfer resources
can generate new concepts.
Generalising the Lanchester `substratum' from a two- to multi-species Lotka-Volterra system
may represent additional actors -- non-combatants,
non-military agencies -- in the operational
environment.  Adding situation awareness to the present three war-fighting functions through Intelligence,
Surveillance and Reconnaissance entities is a next step with a natural spatial
embedding of the model in the `swarmalator' formalism \cite{OKeefe2019}, and may still provide a more compact model than \cite{McLemore2016}.  Recognising that the 
`enemy gets a vote' suggests the
utility of Game Theory here, 
analogous to the treatment of the
Lanchester model in \cite{Morse1951},
with a first step taken in \cite{Demazy2018}. This paper is
thus a significant step in expanding 
equation based approaches
to complex warfare.

\section*{Acknowledgements} 
This work was funded by the Defence Science and Technology (DST) Group of the Department of Defence and conducted under its Modelling Complex Warfighting Strategic Research Initiative. 

\section*{Appendix A: Calculation of reduced steady-state global Lanchester factor}
We begin by substituting the approximation, given by Equation (\ref{approx}), into the defining expression for the Kuramotom-Sakaguchi model (Equation (\ref{eq:BlueKuramoto})) to obtain
\begin{eqnarray}
\begin{split}
\dot{b}_i + \dot{\Theta}_B = \omega_i - \sigma_B f \left(p_B \right) \sum_{j \in {\cal{B}}} L^{B}_{ij}b_j - \zeta_{BR} g\left(p_B \right) \sin \left( \Delta_{BR} - \phi_{BR} \right) d^{BR}_i\\
- \zeta_{BR} g\left(p_B \right) \cos \left( \Delta_{BR} - \phi_{BR} \right)\sum_{j \in {\cal{B}} \cup {\cal{R}} }{\cal{L}}^{BR}_{ij}v_j,\;\; i \in {\cal B},\\
\dot{r}_i + \dot{\Theta}_R = \nu_i - \sigma_R f \left(p_R \right) \sum_{j \in {\cal{R}}} L^{R}_{ij}r_j + \zeta_{RB} g\left(p_R \right) \sin \left( \Delta_{BR} + \phi_{BR} \right) d^{RB}_i\\
+ \zeta_{RB} g\left(p_R \right) \cos \left( \Delta_{BR} + \phi_{RB} \right)\sum_{j \in {\cal{B}} \cup {\cal{R}} }{\cal{L}}^{RB}_{ij}v_j,\;\; i \in {\cal R},
\end{split}
\label{bigapprox}
\end{eqnarray}
where $d^{B}_i$ and $d^R_i$ designate the degree of node $i$ for the Blue and Red network, respectively. Additionally, $L^B$ and $L^R$ are the intra-networks Laplacian matrices of Blue and Red, respectively, given explicitly by,
\begin{equation}
L^B_{ij} = \delta_{ij} d^B_j - {\cal B}_{ij}, \;\; L^R_{ij} = \delta_{ij} d^R_j - {\cal R}_{ij}.
\end{equation}
Correspondingly, $d^{BR}_i$ signifies the number of edges attached to node $i \in {\cal B}$, specifically also connected to nodes from the Red network --- and vice versa for $d^{RB}_i$, \textit{i.e.}
\begin{equation}
d^{BR}_i = \sum_{j \in {\cal B} \cup {\cal R}}{\cal A}^{BR}_{ij}, \;\; d^{RB}_i = \sum_{j \in {\cal B} \cup {\cal R}}{\cal A}^{RB}_{ij}.
\label{localdeg}
\end{equation}
The sum over relevant $i$'s of each of expression in Equation (\ref{localdeg}) give $d_T^{BR},d_T^{RB}$, introduced in Equation (\ref{totaldeg}). Additionally, the quantities ${\cal{L}}^{BR}$ and ${\cal{L}}^{RB}$ are the inter-network graph Laplacian matrices, given by,
\begin{equation}
{\cal{L}}^{BR}_{ij} = \delta_{ij} d^{BR}_j - {\cal A}^{BR}_{ij}, \;\; {\cal{L}}^{RB}_{ij} = \delta_{ij} d^{RB}_j - {\cal A}^{RB}_{ij}.
\end{equation}
We note the importance of the graph Laplacian
in coupled dynamical systems on networks
attributed to the seminal work of \cite{PecCar1998}.
The expression $v_i$ in Equation (\ref{bigapprox}) represents the Blue and Red phase fluctuations, via
\begin{equation}
v_i = \left\{ \begin{array}{ll}
b_i, \;\; i \in {\cal B}\\
r_i, \;\; i \in {\cal R}
\end{array} . \right.
\end{equation}
We employ the intra-network Laplacians in Equation (\ref{bigapprox}) by exploiting their spanning set of orthonormal eigenvectors, given by
\begin{eqnarray}
\begin{split}
e^{(B,\rho)}_i, \;\; \rho \in {\cal B}_E \equiv \{0,1,\dots, |{\cal B}| -1\} , \;\; \sum_{j \in {\cal B}} L^B_{ij} e^{(B,\rho)}_j = \lambda^B_{\rho} e^{(B,\rho)}_i,\\
e^{(R,\rho)}_i, \;\; \rho \in {\cal R}_E \equiv \{0,1,\dots, |{\cal R}| -1\} , \;\; \sum_{j \in {\cal R}} L^R_{ij} e^{(R,\rho)}_j = \lambda^R_{\rho} e^{(R,\rho)}_i.
\end{split}
\end{eqnarray}
We have distinguished the indices applied between the network nodes --- ${\cal B}$ and ${\cal R}$ --- and the eigen-modes --- ${\cal B}_E$ and ${\cal R}_E$. Each eigen-mode for Blue and Red is associated with a corresponding positive semi-definite Laplacian eigenvalue, labelled $\lambda^B$ and $\lambda^R$ respectively. 
One can anticipate \cite{Boll1998} that, as connected graphs, the Blue and Red eigenspaces contain at least one zero-valued Laplacian eigenvalue, labelled specifically as $\lambda^B_0 = \lambda^R_0 = 0$. The eigenvectors corresponding to these eigenvalues, are conveniently given by the (appropriately normalised) vector consisting of entirely unit-valued entries, \textit{i.e.} $e^{(B,0)}_i = 1/\sqrt{|{\cal B}|}$ and $e^{(R,0)}_j = 1/\sqrt{|{\cal R}|}$, $\forall \;\; i \in {\cal B}$ and $j \in {\cal R}$.

We further simplify by setting the cross terms --- $\sum_{j \in {\cal{B}} \cup {\cal{R}} }{\cal{L}}^{BR}_{ij}v_j$ and $\sum_{j \in {\cal{B}} \cup {\cal{R}} }{\cal{L}}^{RB}_{ij}v_j$ --- in Equation (\ref{bigapprox}) to zero. We expand the fluctuations $b_i$ and $r_i$ in the following set of (non-zero) normal eigenmodes
\begin{equation}
b_i = \sum_{\rho \in {\cal B}/\{0\}} e^{(B,\rho)}_i x_{\rho}, \;\; r_i= \sum_{\rho \in {\cal R}/\{0\}} e^{(R,\rho)}_i y_{\rho},
\end{equation}
into Eq.(\ref{bigapprox}) to obtain
\begin{eqnarray}
\begin{split}
\sum_{\rho \in {\cal B}/\{0\}}  e^{(B,\rho)}_i \left[ \dot{x}_{\rho}  + \sigma_B f \left(p_B \right) \lambda^{B}_{\rho}x_{\rho}\right] + \dot{\Theta}_B = \omega_i   - \zeta_{BR} g\left(p_B \right) \sin \left( \Delta_{BR} - \phi_{BR} \right) d^{BR}_i,\\
\sum_{\rho \in {\cal R}/\{0\}} e^{(R,\rho)}_i \left[ \dot{y}_{\rho} + \sigma_R f \left(p_R \right)\lambda^{R}_{\rho} y_{\rho} \right] + \dot{\Theta}_R = \nu_i  + \zeta_{RB} g\left(p_R \right) \sin \left( \Delta_{BR} + \phi_{RB} \right) d^{RB}_i.
\end{split}
\label{bigapprox2}
\end{eqnarray}
We now exploit the orthonomality of the Blue and Red Laplacian eigenvectors in Eq.(\ref{bigapprox2}) to \textit{solve} for each of the dynamic variables. Beginning with the fluctuations $x$ and $y$, we apply the appropriate non-zero eigenvectors ($ e^{(B,\rho_1)}_i$ for Blue and $ e^{(R,\rho_2)}_i$ for Red) to both sides of Eq.(\ref{bigapprox2}) and sum both the Blue and Red expressions over all nodes to obtain
\begin{eqnarray}
\begin{split}
\dot{x}_{\rho_1} = \sum_{i \in {\cal B}} \left[ \omega_i - \zeta_{BR} g\left(p_B \right) \sin \left( \Delta_{BR} - \phi_{BR} \right) d^{BR}_i \right] e^{(B,\rho_1)}_i - \sigma_B f(p_B) \lambda^B_{\rho_1} x_{\rho_1} , \;\; \rho_1 \in {\cal B}_E / 0\\
\dot{y}_{\rho_2} = \sum_{i \in {\cal R}} \left[ \nu_i + \zeta_{RB} g\left(p_R \right) \sin \left( \Delta_{BR} + \phi_{RB} \right) d^{RB}_i \right] e^{(R,\rho_2)}_i - \sigma_R f(p_R) \lambda^R_{\rho_2} y_{\rho_2} , \;\; \rho_2 \in {\cal R}_E /0
\end{split}
\label{bigapprox3}
\end{eqnarray}
where we have applied the relations,
\begin{equation}
\dot{\Theta}_B = \sqrt{|{\cal B}|}  \dot{\Theta}_B e^{(B,0)}_i, \;\; \dot{\Theta}_R = \sqrt{|{\cal R}|}  \dot{\Theta}_R e^{(R,0)}_i.
\end{equation}
Equation (\ref{bigapprox3}) gives a convenient linear dynamic expression for each of the fluctuation variables, $x$ and $y$. We note however that solving for either $x$ or $y$ requires the dynamic expression for the difference of the global phases, $\Delta_{BR}$, in addition to $f$ and $g$, coming from the Lanchester dynamics. 

Focusing instead on the global phases $\Theta_B$ and $\Theta_R$, as opposed to the fluctuations, we apply the appropriate zero eigenmodes ($ e^{(B,0)}_i$ for Blue and $ e^{(R,0)}_i$ for Red) to both sides of Eq.(\ref{bigapprox2}), and sum over the appropriate (Blue or Red) network nodes to obtain
\begin{equation}
\dot{\Theta}_B = \bar{\omega}_B - \frac{\zeta_{BR}}{|{\cal B}|} g(p_B) \sin (\Delta_{BR} - \phi_{BR}) d^{BR}_T,\;\; 
\dot{\Theta}_R = \bar{\omega}_R + \frac{\zeta_{RB}}{|{\cal R}|} g(p_R) \sin (\Delta_{BR} + \phi_{RB}) d^{RB}_T,
\label{THETA2}
\end{equation} 
where $\bar{\omega}_B, \bar{\omega}_R$ were introduced
earlier as the mean of the Blue and Red natural frequencies.  Hence, taking the difference of the time derivatives of the two global frequencies in Equation (\ref{THETA2}) awards us with the defining expression for $\dot{\Delta}_{BR}$. Combining this with the appropriately approximated global Lanchester expressions (Equation \ref{fullLan}), \textit{i.e.} $O^{global}_{B} \approx O^{global}_{R} \approx 1$, our dimensionally reduced system is given by Equation (\ref{totalapprox}), and restated below for convenience
\begin{eqnarray} 
\begin{split}
\dot{\Delta}_{BR} = \bar{\omega}_B-\bar{\omega}_R -\left[ \frac{\zeta_{BR}}{|{\cal B}|} g(p_B) \sin (\Delta_{BR} - \phi_{BR}) d^{BR}_T +\frac{\zeta_{RB}}{|{\cal R}|} g(p_R) \sin (\Delta_{BR} + \phi_{RB}) d^{RB}_T \right],\\
\dot{p}_B = -\kappa_{RB} \frac{1-\sin(\Delta_{BR}) }{2}p_R\mathcal{H}(p_B), \;\;
\dot{p}_R = -\kappa_{BR} \frac{1+\sin(\Delta_{BR}) }{2}p_B\mathcal{H}(p_R).
\end{split}
\label{totalapprox22}
\end{eqnarray}
Noting that, in contrast to the corresponding time derivative expressions for the fluctuations in Equation (\ref{bigapprox3}), the time derivative expression for the difference in global phases given in Equation (\ref{totalapprox}) is independent of the fluctuations $x$ and $y$, in addition to the Lanchester based feedback term $f$; thus in the regimes where $O_B \approx O_R \approx 1$, the presence of $f$ is immaterial, explaining our observation that different $f$ values led to contours similar to figure \ref{fig:CombinationPlot}.

\section*{Appendix B: Detailed
exploration of use-case 2}
To explore the different effects in 
a scenario as simple as use-case 2 detailed in section \ref{fightfish}, in Fig.\ref{fig:FFLabelledPhasePlot}
we set an intermediate inter-network coupling $\zeta_{RB}=\zeta_{BR}=1.25$
and traverse a particular path through the landscape as follows: 

\begin{figure}[h]
	\centering
	\includegraphics[width=0.5\linewidth]{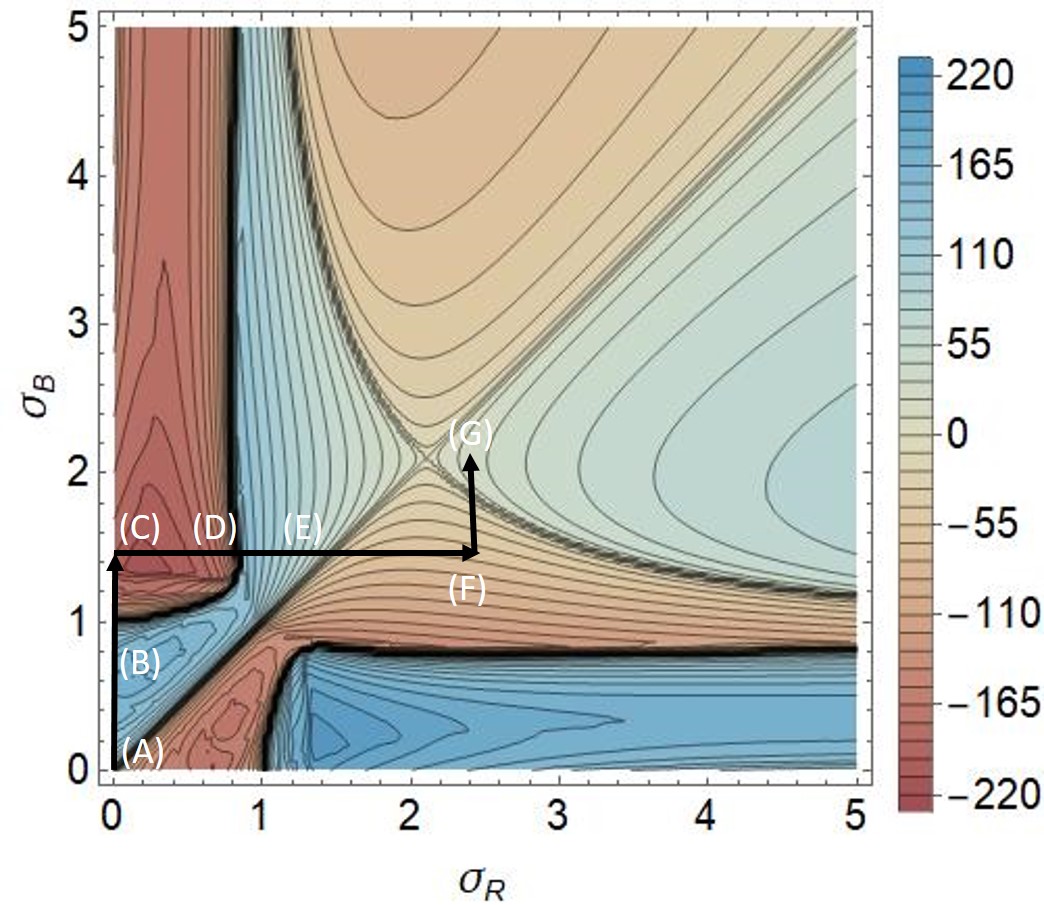}
	\caption{A phase plot of combat outcomes with a labelled path. $\zeta_{BR}=\zeta_{RB}=1.25$.}
	\label{fig:FFLabelledPhasePlot}
\end{figure}

\begin{itemize}
	\item Point (A) with $(\sigma_R,\sigma_B) = (0,0)$: Starting at the origin, this is the line of neutral advantage $\sigma_B=\sigma_R$, and thus neither side is able to achieve supremacy. With zero intra-coupling, engagement nodes traverse the decision cycle at a higher frequency than their counterparts due to the frustrations causing a forcing effect.
	\item Path (A-B): As Blue's internal coupling increases, the engagement and manoeuvre nodes begin to influence each other. Nevertheless, the forcing due to frustration on engagement nodes is still too large for internal coupling to completely overcome, resulting in periodic behaviour of engagement nodes decoupling and re-synchronising with the bulk.  Notably, as Blue's internal coupling increases, the engagement nodes spend more time per cycle in the period before decoupling than afterwards.
	\item Point (B) with $(\sigma_R,\sigma_B) = (0,0.6)$: In this region Blue obtains an optimal balance between receiving resources fast enough to achieve a numerical advantage and not allowing Red to attain an overwhelming decision (phase) advantage.
	\item Path (B-C): Further increases to Blue's internal coupling continue to afford Red a larger phase advantage whilst the manoeuvre benefits begin to plateau. Once $\sigma_B > 1$, Red begins to achieve victory in engagements.
	\item Point (C) with $(\sigma_R,\sigma_B) = (0,1.5)$: $\sigma_B$ is large enough that the Blue network no longer decouples, and Red only receives resource and a higher local order parameter sporadically when combat nodes happen to coincide. However, this is more than sufficient to overcome the initial boost in Blue's superiority in resource transfer. Additionally, given Blue is only synchronised to a small extent, the local order parameter is often suboptimal as the non-engagement nodes frequency synchronise in a lagged position. Conversely, although Red spends a significant time with an even lower local order parameter between engagement and manoeuvre nodes, this nevertheless results in devastating losses for Blue.
	\item Path (C-D): A slight increase in Red internal coupling allows the Red non-combat nodes to group up, facilitating a faster Red manoeuvre when they happen to phase align with Red combat nodes. With the most optimal region of Red combat results being the domain around point (D).
	\item Path (D-E): As $\sigma_R$ increases, Red begins to sacrifice its favourable engagement outcomes, with the critical transition occurring around $\sigma_R\approx0.8$, when its engagement and manoeuvre nodes no longer periodically decouple. Past this value, both Blue and Red reach steady group frequencies rather than periodic decoupling behaviour. Notably, the increased Blue manoeuvre and local order parameter afforded by its larger internal coupling dominates the now reduced Red phase advantage.
	\item Path (E-F): Once $\sigma_R$ crosses the neutral phase advantage line $\sigma_R=\sigma_B$, its larger intra-network coupling allows it to achieve victory over Blue. At this point both networks are receiving diminishing return on manoeuvre benefits and so the optimisation is between local order parameter between internal nodes and phase advantage over the opposition's engagement nodes.
	\item Path (F-G): Increasing $\sigma_B$ to Point (G) $(\sigma_R,\sigma_B) = (2.5,2.1)$ allows Blue to gain a phase advantage over Red which provides a slight advantage that grows over the course of the engagement.
\end{itemize}

In summary, intra and inter-network coupling work to divide the landscape into a number of regions according to whether internal coupling is large enough to overcome the tension between engagement and manoeuvre nodes imposed by the value of $\zeta$.

\end{document}